%% file: m3-kato.tex
\input m3-macs

\pageno=165
 
\tinfos {}.{}.165-195

\let\HHsav\HH
\def\HH{\vglue .1cm \HHsav}
\let\HHHsav\HHH
\def\HHH{\vglue .1cm \HHHsav}

\SetTFLinebox{\gtp }
\SetFLinebox{\gtv3 }
\SetHLinebox{\issn}

\vglue .3cm

\H Existence theorem for higher local fields

Kazuya Kato

\SetAuthorHead{K. Kato}
\SetTitleHead{Existence theorem for higher local fields}

\vglue .3cm

\HH 0. Introduction

\vglue .2cm
A field $K$ is called an $n$-dimensional local field if there is
a sequence of fields $k_n,\dots, k_0$ satisfying the following conditions:
$k_0$ is a finite field, $k_i$ is a complete discrete valuation field
with residue field $k_{i-1}$ for $i=1,\dots,n$,
and $k_n=K$.

In \cite{9} we defined a canonical homomorphism from
the $n$th Milnor group $K_n(K)$ (cf.\ \cite{14}) of an $n$-dimensional local field $K$
to the Galois group $\Gal(K^{\ab}/K)$ of the maximal abelian extension of $K$
and generalized the familiar results of the usual local class field theory
to the case of arbitrary dimension except the ``existence theorem''.

An essential difficulty with the existence theorem
lies in the fact that $K$ (resp.\ the multiplicative group $K^*$)
has no appropriate topology in the case where $n\ge 2$ (resp.\ $n\ge 3$)
which would be compatible with the ring (resp.\ group) structure
and which would take the topologies of the residue fields
into account.
Thus we  abandon the familiar tool ``topology''
and define the openness of subgroups and the continuity of maps
from a new point of view.

In the following main theorems the words ``open'' and ``continuous''
are not used in the topological sense.
They are explained below.

\th Theorem 1

Let $K$ be an $n$-dimensional local field.
Then the correspondence $$L\to N_{L/K}K_n(L)$$ is a bijection from the set
of all finite abelian extensions of $K$ to the set of all open subgroups of $K_n(K)$ of finite index.
\endth

This existence theorem is essentially contained in the following theorem
which expresses certain Galois cohomology groups of $K$ (for example
the Brauer group of $K$) by using the Milnor $K$-group of $K$.
For a field $k$ we define the group $H^r(k)$ ($r\ge0$)
as follows (cf.\ \cite{9,\S 3.1}).
In the case where $\chr(k)=0$ let 
$$H^r(k)=\inlim H^r(k,\mu_m^{\otimes(r-1)})$$
(the Galois cohomology).
In the case where $\chr(k)=p>0$ let
$$H^r(k)=\inlim H^r(k, \mu_m^{\otimes(r-1)})+\inlim H_{p^i}^r(k).$$
Here in each case $m$ runs over all integers invertible
in $k$,
$\mu_m$ denotes the group of all $m$th roots of 1 in the separable closure 
$k^{\sep}$ of $k$,
and $\mu_m^{\otimes(r-1)}$ denotes its $(r-1)$th tensor power as a $\Bbb Z/m$-module 
on which $\Gal(k^{\sep}/k)$ acts in the natural way.
In the case where $\chr(k)=p>0$ we denote
by $H_{p^i}^r(k)$ the cokernel of
$$F-1\colon C_i^{r-1}(k)\to C_i^{r-1}(k)/\{C_i^{r-2}(k),T\}$$ 
where $C_i^{\cdot}$ is the group defined in \cite{3, Ch.II,\S 7}
(see also Milne \cite{13, \S3}).
For example, $H^1(k)$ is isomorphic to the group of all continuous characters of the compact abelian group
$\Gal(k^{\ab}/k)$
and $H^2(k)$ is isomorphic to the Brauer group of $k$.

\th Theorem 2

Let $K$ be as in Theorem 1.
Then $H^r(K)$ vanishes for $r>n+1$
and is isomorphic to the group of all continuous characters of finite order
of $K_{n+1-r}(K)$ in  the case where $0\le r\le n+1$.
\endth

We shall explain the contents of each section.

For a category $\Cal C$ the category of pro-objects $\pro(\Cal C)$
and the category of ind-objects $\ind(\Cal C)$ 
are defined as in Deligne \cite{5}.
Let $\Cal F_0$ be the category of finite sets,
and let $\Cal F_1$, $\Cal F_2$, \dots be the categories defined
by $\Cal F_{n+1}=\ind(\pro(\Cal F_n))$. Let $\Cal F_\infty=\cup_n \Cal F_n$.
In section 1 we shall show that $n$-dimensional local fields
can be viewed as ring objects of $\Cal F_n$.
More precisely we shall define a ring object
$\U K$ of $\Cal F_n$ corresponding to an $n$-dimensional local field $K$
such that $K$ is identified with the ring $[e,\U K]_{\Cal F_\infty}$ of morphisms from
the one-point set $e$ (an object of $\Cal F_0$) to $\U K$,
and a group object $\U{K^*}$ such that $K^*$ is identified with
$[e,\U {K^*}]_{\Cal F_\infty}$.
We call a subgroup $N$ of $K_q(K)$ open if and only if the map
$$K^*\times \dots\times K^* \to K_q(K)/N,\quad (x_1,\dots,x_q)\mapsto \{x_1,\dots,x_q\}\mod N$$ comes from a morphism 
$\U{K^*}\times \dots\times \U{K^*} \to K_q(K)/N$
of $\Cal F_\infty$
where $K_q(K)/N$ is viewed as an object of $\ind(\Cal F_0)\subset \Cal F_1$.
We call a homomorphism $\varphi\colon K_q(K)\to \Bbb Q/\Bbb Z$
a continuous character if and only if
the induced map
$$K^*\times \dots\times K^* \to \Bbb Q/\Bbb Z,\quad
(x_1,\dots,x_q)\mapsto  \varphi(\{x_1,\dots,x_q\})$$
comes from a morphism of $\Cal F_\infty$
where $\Bbb Q/\Bbb Z $ is viewed as an object of $\ind(\Cal F_0)$.
In each case such a morphism of $\Cal F_\infty$  
is unique if it exists (cf.\ Lemma 1 of section 1).

In section 2 we shall generalize the self-duality of the additive group of 
a one-dimensional local field in the sense of Pontryagin
to arbitrary dimension.

Section 3 is a preliminary one for section 4.
There we shall prove some ring-theoretic properties of
$[X,\U K]_{\Cal F_\infty}$ for objects $X$ of $\Cal F_\infty$.

In section 4 we shall treat the norm groups of cohomological objects.
For a field $k$ denote by $\Cal E(k)$ the category of all finite extensions
of $k$ in a fixed algebraic closure of $k$
with the inclusion maps as morphisms.
Let $H$ be a functor from $\Cal E(k)$ to the category $\text{\rm Ab}$ of all abelian groups
such that
$\inlim_{k'\in \Cal E(k)} \, H(k')=0$.
For $w_1,\dots,w_g\in H(k)$ define the $K_q$-norm group
$N_q(w_1,\dots,w_g)$ as the subgroup of $K_q(k)$
generated by the subgroups
$N_{k'/k} K_q(k')$ where $k'$ runs over all fields in $\Cal E(k)$ such that
$\{w_1,\dots,w_g\} \in \kr(H(k)\to H(k'))$ and 
where $N_{k'/k}$ denotes the canonical norm homomorphism of the Milnor  
$K$-groups
(Bass and Tate \cite{2, \S5} and \cite{9, \S1.7}).
For example, if $H=H^1$ and
$\chi_1,\dots,\chi_g\in H^1(k)$ then
$N_q(\chi_1,\dots,\chi_g)$ is nothing but
$N_{k'/k} K_q(k')$ where $k'$ is the finite abelian extension
of $k$ corresponding to 
$\cap_i \kr(\chi_i\colon \Gal(k^{\ab}/k)\to \Bbb Q/\Bbb Z)$.
If $H=H^2$ and $w\in H^2(k)$ then
$N_1(w)$ is the image of the reduced norm map $A^*\to k^*$
where $A$ is a central simple algebra over $k$ corresponding to $w$.

As it is well known for a one-dimensional local field $k$
the group $N_1(\chi_1,\dots,\chi_g)$ is an open subgroup of $k^*$
of finite index for any $\chi_1,\dots,\chi_g\in H^1(k)$
and the group $N_1(w)=k^*$ for any $w\in H^2(k)$.
We generalize these facts as follows.

\th Theorem 3

Let $K$ be an $n$-dimensional local field and let $r\ge 1$.

\Roster 
\Item{(1)} Let $w_1,\dots, w_g \in H^r(K)$.
Then the norm group
$N_{n+1-r}(w_1,\dots,w_g)$ is an open subgroup of
$K_{n+1-r}(K)$ of finite index.

\Item{(2)} Let $M$ be a discrete torsion abelian group endowed with a continuous
action of $\Gal(K^{\sep}/K)$.
Let $H$ be the Galois cohomology functor
$H^r(\,\,\,,M)$.
Then for every $w\in H^r(K,M)$ the group
$N_{n+1-r}(w)$ is an open subgroup of $K_{n+1-r}(K)$ of finite index.
\endRoster
\endth

Let $k$ be a field and let $q,r\ge 0$.
We define a condition $(N_q^r,k)$ as follows: 
for every $k'\in \Cal E(k)$ and every discrete
torsion abelian group $M$ endowed with a continuous action of $\Gal({k'}^{\sep}/k')$
$$N_q(w_1,\dots,w_g)=K_q(k')$$
for every $i>r$,
$w_1,\dots,w_g\in H^i(k')$, $w_1,\dots, w_g\in H^i(k',M)$,
and in addition $|k:k^p|\le p^{q+r}$ in  the case where $\chr(k)=p>0$.

For example, if $k$ is a perfect field then the condition $(N_0^r,k)$
is equivalent to $\cd(k)\le r$ where $\cd$ denotes the cohomological dimension 
(Serre \cite{16}).

\th Proposition 1

Let $K$ be a complete discrete valuation field
with residue field $k$.
Let $q\ge 1$ and $r\ge 0$.
Then the two conditions
$(N_q^r,K)$ and $(N_{q-1}^r,k)$
are equivalent.
\endth

On the other hand by \cite{11} 
the conditions $(N_0^r,K)$ and $(N_0^{r-1}, k)$
are equivalent for any $r\ge 1$.
By induction on $n$ we obtain

\th Corollary

Let $K$ be an $n$-dimensional local field.
Then the condition $(N_q^r,K)$ holds
if and only if $q+r\ge n+1$.
\endth

We conjecture that
if $q+r=q'+r'$ then the two conditions
$(N_q^r,k)$ and $(N_{q'}^{r'},k)$
are equivalent for any field $k$.

\smallskip

Finally in section 5 we shall prove Theorem 2.
Then Theorem 1 will be a corollary of Theorem 2  for $r=1$
and of \cite{9,\S3, Theorem 1} which claims
that the canonical homomorphism
$$K_n(K)\to \Gal(K^{\ab}/K)$$ 
induces an isomorphism
$K_n(K)/N_{L/K}K_n(L)\iss \Gal (L/K)$
for each finite abelian extension $L$ of $K$.

I would like to thank Shuji Saito for helpful discussions
and for the stimulation given by his research in this area
(e.g. his duality theorem of Galois cohomology groups with locally compact
topologies for two-dimensional local fields).

\HHH Table of contents

\phantom{}\par 

1. Definition of the continuity for higher local fields.

2. Additive duality.

3. Properties of the ring of $\U K$-valued morphisms.

4. Norm groups.

5. Proof of Theorem 2.

\HHH Notation

\phantom{}\par

We follow the notation in the beginning of this volume.
References to sections in this text mean
references to sections of this work and not of the whole volume. 

All fields and rings in this paper are assumed to be commutative.

Denote by $\text{\rm Sets, Ab, Rings}$ the categories of sets, of abelian groups and of rings respectively.

If $\Cal C$ is a category and $X,Y$ are objects of $\Cal C$
then $[X,Y]_{\Cal C}$ (or simply $[X,Y]$)
denotes the set of morphisms $X\to Y$.

\bigskip

\HH 1. Definition of the continuity for higher local fields 

\HHH 1.1. Ring objects of a category corresponding to rings

\phantom{}\smallskip\par

For a category $\Cal C$ let $\Cal C^{\circ}$ be the dual category
of $\Cal C$.
If $\Cal C$ has a final object we always denote it by $e$.
Then, if $\theta\colon X\to Y$ is a morphism of $\Cal C$,
$[e,\theta]$ denotes the induced map $[e,X]\to [e,Y]$.

In this subsection we prove the following 

\th Proposition 2

Let $\Cal C$ be a category with a final object $e$ in which
the product of any two objects exists.
Let $\U R$ be a ring object of $\Cal C$ such that
for a prime $p$ the morphism
$\U R\to \U R$, $x\mapsto px$ is the zero morphism,
and via the morphism $\U R\to \U R$, $x\mapsto x^p$ the latter $\U R$ is a free module
of finite rank over the former $\U R$.
Let $R=[e,\U R]$, and let
$A$ be a ring with a nilpotent ideal $I$ such that $R=A/I$ and such that
$I^i/I^{i+1}$ is a free $R$-module of finite rank for any $i$.

Then:
\Roster 
\Item{(1)} There exists a ring object $\U A$ of $\Cal C$ equipped
with a ring isomorphism $j\colon A\iss [e,\U A]$
and with a homomorphism of ring objects $\theta\colon \U A\to \U R$ having
the following properties:
\ItemItem{(a)} $[e,\theta]\circ j\colon A\to R$ coincides with
the canonical projection.
\ItemItem{(b)} For any object $X$ of $\Cal C$, $[X,\U A]$ is a formally etale ring over $A$
in the sense of Grothendieck \cite{7, Ch. 0 \S19},
and $\theta$ induces an isomorphism
$$[X,\U A]/I[X,\U A]\simeq [X,\U R].$$ 
 \Item{(2)} The above triple $(\U A, j,\theta)$ is unique
in the following sense.
If $(\U A', j',\theta')$ is another triple satisfying the same condition in (1),
then there exists a unique isomorphism of ring objects
$\psi\colon \U A\iss \U A'$ such that $[e,\psi]\circ j = j'$
and $\theta=\theta'\circ \psi$.

\Item{(3)} The object $\U A$ is isomorphic (if one forgets the ring-object
structure)
to the product of finitely many copies of $\U R$.

\Item{(4)} If $\Cal C$ has finite inverse limits, the above
assertions (1) and (2)
are valid if  conditions ``free module of finite rank'' on $\U R$
and $I^i/I^{i+1}$ are replaced by  conditions ``direct summand of a free
module of finite rank''.
\endRoster 
\endth

\eg Example

Let $R$ be a non-discrete locally compact field and $A$ a local ring
of finite length with residue field $R$.
Then in  the case where $\chr(R)>0$ Proposition 2
shows that there exists a canonical topology on $A$ compatible with the ring
structure
 such that $A$ is homeomorphic to the product of finitely many copies of $R$.
On the other hand, in the  case where $\chr(R)=0$ it is impossible in general
to define canonically such a topology on $A$.
Of course, by taking a section $s\colon R\to A$ (as rings),
$A$ as a vector space over $s(R)$ has the vector space topology,
but this topology depends on the choice of $s$ in general.
This reflects the fact that in  the case of $\chr(R)=0$ the ring of $R$-valued continuous functions on a topological space is not in general 
formally smooth over $R$  contrary to the case of $\chr(R)>0$.
\endeg

\pf Proof of Proposition 2

Let $X$ be an object of $\Cal C$; put
$R_X=[X,\U R]$.
The assumptions on $\U R$ show that the homomorphism
$$R^{(p)}\otimes_R R_X\to R_X,\quad x\otimes y \mapsto xy^p$$
is bijective, where
$R^{(p)}=R$ as a ring and the structure homomorphism
$R\to R^{(p)}$ is $x\mapsto x^p$.
Hence by \cite{10, \S1 Lemma 1} there exists a formally
etale ring $A_X$ over $A$ with a ring isomorphism 
$\theta_X\colon A_X/IA_X\simeq R_X$.
The property ``formally etale'' shows that the correspondence $X\to A_X$ is a functor
$\Cal C^\circ \to \text{\rm Rings}$, and that the system $\theta_X$ forms a morphism
of functors.
More explicitly, let $n$ and $r$ be sufficiently large integers,
let $W_n(R)$ be the ring of $p$-Witt vectors over $R$ of length $n$,
and let $\varphi\colon W_n(R)\to A$ be the homomorphism
$$(x_0,x_1,\dots)\mapsto \sum_{i=0}^r p^i\wt{x_i}^{p^{r-i}}$$
where $\wt{x_i}$ is a representative of $x_i\in R$ in $A$.
Then $A_X$ is defined as 
the tensor product
$$W_n(R_X)\otimes_{W_n(R)} A$$ induced by $\varphi$.
Since $\text{\rm Tor\,}_1^{W_n(R)} (W_n(R_X), R)=0$ we have
$$\text{\rm Tor\,}_1^{W_n(R)} (W_n(R_X), A/I^i)=0$$ for every $i$.
This proves that the canonical homomorphism
$$I^i/I^{i+1}\otimes_{R} R_X \to I^iA_X/I^{i+1}A_X$$
is bijective for every $i$.
Hence each functor $X\to I^iA_X/I^{i+1}A_X$ is representable by a finite
product of copies of $\U R$,
and it follows immediately that the functor $A_X$ is represented by the product
of finitely many  copies of $\U R$.
\qed\endpf 

\HHH 1.2. $n$-dimensional local fields as objects of $\Cal F_n$

\phantom{}\smallskip\par

Let $K$ be an $n$-dimensional local field.
In this subsection we define a ring object $\U K$
and a group object $\U{K^*}$ by induction on $n$.

Let $k_0,\dots,k_n=K$ be as in the introduction.
For each $i$ such that $\chr(k_{i-1})=0$ (if such an $i$ exists)
choose a ring morphism $s_i\colon k_{i-1}\to \Cal O_{k_i}$
such that the composite
$k_{i-1}\to \Cal O_{k_i}\to \Cal O_{k_i}/\Cal M_{k_i}$ is the indentity map.
Assume $n\ge 1$ and let
$\U{k_{n-1}}$ be the  ring object of $\Cal F_{n-1}$ corresponding to $k_{n-1}$
by induction on $n$.

If $\chr(k_{n-1})=p>0$, the construction of $\U K$ below will show
by induction on $n$ that the assumptions of Proposition 2 are satisfied when one takes
$\Cal F_{n-1}$, $\U{k_{n-1}}$, $k_{n-1}$ and $\Cal O_K/\Cal M_K^r$ ($r\ge 1$)
as $\Cal C$, $\U R$, $R$ and $A$.
Hence we obtain a ring object $\U{\Cal O_K/\Cal M_K^r}$ of
$\Cal F_{n-1}$.
We identify $\Cal O_K/\Cal M_K^r$ with
$[e,\U{\Cal O_K/\Cal M_K^r}]$ via the isomorphism $j$ of Proposition 2.

If $\chr(k_{n-1})=0$, let $\U{\Cal O_K/\Cal M_K^r}$
be the ring object of $\Cal F_{n-1}$ which represents
the functor
$$\Cal F_{n-1}^{\circ}\to \text{\rm Rings},\quad X\mapsto \Cal O_K/\Cal M_K^r\otimes_{k_{n-1}}[X,\U{k_{n-1}}],$$
where $\Cal O_K/\Cal M_K^r$ is viewed as a ring over $k_{n-1}$ via $s_{n-1}$.

In each case let $\U{\Cal O_K}$ be the object
"$\varprojlim$"$\U{\Cal O_K/\Cal M_K^r}$ of $\pro(\Cal F_{n-1})$.
We define $\U K$ as the ring object of $\Cal F_n$
which corresponds to the functor
$$\pro(\Cal F_{n-1})^{\circ}\to \text{\rm Rings},
\quad X\mapsto K\otimes_{\Cal O_K}[X,\U{\Cal O_K}].$$

Thus, $\U K$ is defined canonically in  the case of $\chr(k_{n-1})>0$,
and it depends (and doesn't depend) on the choices of $s_i$
in  the case of $\chr(k_{n-1})=0$ in the following sense.
Assume that another choice of sections $s_i'$ yields $\U{k_i}'$
and $\U{K}'$.
Then there exists an isomorphism
of ring objects $\U K\iss \U{K}'$ which induces
$\U{k_i}\iss \U{k_i}'$ for each $i$.
But in general there is no isomorphism of ring objects
$\psi\colon \U K\to \U{K}'$ such that
$[e,\psi]\colon K\to K$ is the indentity map.

Now let $\U{K^*}$ be the object of $\Cal F_n$ which represents
the functor
$$\Cal F_n^{\circ}\to \text{\rm Sets},\quad X\mapsto [X,\U K]^*.$$
This functor is representable because
$\Cal F_n$ has  finite inverse limits as can be shown by induction on $n$.

\df Definition 1

We define fine (resp.\ cofine) objects of $\Cal F_n$ by induction on $n$.
All objects in $\Cal F_0$ are called fine (resp.\ cofine) objects of $\Cal F_0$.
An object of $\Cal F_n$ ($n\ge 1$) is called a fine (resp.\ cofine) object
of $\Cal F_n$ if and only if it is expressed as
$X=$"$\inlim$"$X_\lambda$ for some objects $X_\lambda$
of $\pro(\Cal F_{n-1})$ and each $X_\lambda$ is expressed as
$X_\lambda=$"$\varprojlim$"$X_{\lambda \mu}$ for some objects $X_{\lambda\mu}$
of $\Cal F_{n-1}$ satisfying the condition that all $X_{\lambda\mu}$ are fine
(resp.\ cofine) objects of $\Cal F_{n-1}$
and the maps $[e,X_\lambda]\to  [e,X_{\lambda\mu}]$ are surjective for all
$\lambda,\mu$ (resp.\ the maps $[e,X_\lambda]\to [e,X]$ are injective for all $\lambda$).
\enddf

Recall that if $i\le j$ then $\Cal F_i$ is a full subcategory of $\Cal F_j$.
Thus each $\Cal F_i$ is a full subcategory of $\Cal F_\infty=\cup_i\Cal F_i$.

\th Lemma 1

\Roster 
\Item{(1)} Let $K$ be an $n$-dimensional local field.
Then an object of $\Cal F_n$ of the form
$$\U K\times \dots \U K\times\U{K^*}\times\dots\times \U{K^*}$$ 
is a fine and cofine
object of $\Cal F_n$.
Every set $S$ viewed as an object of $\ind(\Cal F_0)$ 
is a fine and cofine object of $\Cal F_1$.

\Item{(2)} Let $X$ and $Y$ be objects of $\Cal F_\infty$,
and assume that $X$ is a fine object of $\Cal F_n$
for some $n$ and $Y$ is a cofine object of $\Cal F_m$
for some $m$.
Then two morphisms $\theta,\theta'\colon X\to Y$ coincide
if $[e,\theta]=[e,\theta']$.
\endRoster 
\endth

As explained in 1.1 the definition of the object $\U K$
depends on the sections $s_i\colon k_{i-1}\to \Cal O_{k_i}$
chosen for each $i$ such that
$\chr(k_{i-1})=0$.
Still we have the following: 

\th Lemma 2

\Roster 
\Item{(1)} Let $N$ be a subgroup of $K_q(K)$ of finite index.
Then openness of $N$ doesn't depend on the choice of sections $s_i$.

\Item{(2)} Let $\varphi\colon K_q(K)\to \Bbb Q/\Bbb Z$ be a homomorphism
of finite order.
Then the continuity of $\chi$ doesn't depend on the choice of sections $s_i$.
\endRoster
\endth

The exact meaning of Theorems 1,2,3 is now clear.

\HH 2. Additive duality

\HHH 2.1. Category of locally compact objects

\phantom{}\smallskip\par

If $\Cal C$ is the category of finite abelian groups,
let $\wt{\Cal C}$ be the category of topological abelian groups $G$ which possess a totally disconnected open compact subgroup $H$
such that $G/H$ is a torsion group.
If $\Cal C$ is the category of finite dimensional vector spaces over
a fixed (discrete) field $k$, let $\wt{\Cal C}$ be the category
of locally linearly compact vector spaces over $k$ (cf.\ Lefschetz \cite{12}).
In both cases the canonical self-duality
of $\wt{\Cal C}$ is well known.
These two examples are special cases of the following general construction.

\df Definition 2

For a category $\Cal C$ define a full subcategory
$\wt{\Cal C}$ of $\ind(\pro(\Cal C))$
as follows.
An object $X$ of $\ind(\pro(\Cal C))$ belongs to $\wt{\Cal C}$
if and only if
it is expressed in the form
"$\inlim$"${}_{j\in J}$"$\varprojlim$"${}_{i\in I} X(i,j)$ for some directly ordered
sets $I$ and $J$ viewed as small categories in the usual way
and for some functor $X\colon I^{\circ}\times J\to \Cal C$
satisfying the following conditions.
\Roster 
\Item{(i)} If $i,i'\in I$, $i\le i'$
then the morphism $X(i',j)\to X(i,j)$ is surjective
for every $j\in J$.
If $j,j'\in J$, $j\le j'$
then the morphism $X(i,j)\to X(i,j')$ is injective
for every $i\in I$.

\Item{(ii)} If $i,i'\in I$, $i\le i'$
and $j,j'\in J$, $j\le j'$ then the square
$$
\CD
X(i',j) @>>> X(i',j')\\
@VVV @VVV \\
X(i,j)@>>> X(i,j')
\endCD
$$
is cartesian and cocartesian.
\endRoster
\enddf

It is not difficult to prove that $\wt{\Cal C}$ is equivalent to the
full subcategory of $\pro(\ind(\Cal C))$ (as well as $\ind(\pro(\Cal C))$)
consisting of all objects which are expressed in the form
 "$\varprojlim$"${}_{i\in I}$"$\inlim$"${}_{j\in J} X(i,j)$ 
for some triple $(I,J,X)$ satisfying the same conditions as above.
In this equivalence the object 

\noindent "$\inlim$"${}_{j\in J}$"$\varprojlim$"${}_{i\in I} X(i,j)$ corresponds to 
 "$\varprojlim$"${}_{i\in I}$"$\inlim$"${}_{j\in J} X(i,j)$. 

\df Definition 3

Let $\Cal A_0$ be the category of finite abelian groups,
and let $\Cal A_1,\Cal A_2,\dots$ be the categories defined
as $\Cal A_{n+1}=\wt{\Cal A_n}$.
\enddf

It is easy to check by induction on $n$
that $\Cal A_n$ is a full subcategory of the category
$\Cal F_n^{\ab}$ of all abelian group objects of $\Cal F_n$
with additive morphisms.

\HHH 2.2. Pontryagin duality

\phantom{}\smallskip\par

The category $\Cal A_0$ is equivalent to its dual via the functor
$$D_0\colon \Cal A_0^{\circ}\iss \Cal A_0,
\quad X\mapsto \Hom(X,\Bbb Q/\Bbb Z).$$
By induction on $n$ we get an equivalence
$$D_n\colon \Cal A_n^{\circ}\iss \Cal A_n,
\quad \Cal A_n^{\circ}=(\wt{\Cal A_{n-1}})^{\circ}=
\wt{\Cal A_{n-1}^{\circ}}@> D_{n-1}>> \wt{\Cal A_{n-1}}=\Cal A_n$$
where we use $(\wt{\Cal C})^{\circ}=
\wt{\Cal C^{\circ}}$.
As in  the case of $\Cal F_n$
each $\Cal A_n$ is a full subcategory of $\Cal A_\infty=\cup_n \Cal A_n$.
The functors $D_n$ induce an equivalence
$$D\colon \Cal A_\infty^\circ\iss \Cal A_\infty$$
such that $D\circ D$ coincides with the indentity functor.

\th Lemma 3

View $\Bbb Q/\Bbb Z$ as an object of $\ind(\Cal A_0)\subset \Cal A_\infty\subset \Cal F_\infty^{\ab}$.
Then:
\Roster 
\Item{(1)} For every object $X$ of $\Cal A_\infty$
$$[X,\Bbb Q/\Bbb Z]_{\Cal A_\infty}\simeq [e,D(X)]_{\Cal F_\infty}.$$

\Item{(2)}  For all objects $X,Y$ of $\Cal A_\infty$
$[X,D(Y)]_{\Cal A_\infty}$ is canonically isomorphic to
the group of biadditive morphisms $X\times Y\to \Bbb Q/\Bbb Z$ in
$\Cal F_\infty$.
\endRoster
\endth
\pf Proof

The isomorphism of (1) is given by
$$[X,\Bbb Q/\Bbb Z]_{\Cal A_\infty}\simeq
[D(\Bbb Q/\Bbb Z),D(X)]_{\Cal A_\infty}=
[\wh{\Bbb Z},D(X)]_{\Cal A_\infty}\iss [e,D(X)]_{\Cal F_\infty}$$
($\wh{\Bbb Z}$ is the totally disconnected compact abelian group
$\varprojlim_{n>0} \Bbb Z/n$ and the last arrow
is the evaluation at $1\in \wh{\Bbb Z}$).
The isomorphism of (2) is induced by the canonical
biadditive morphism $D(Y)\times Y\to \Bbb Q/\Bbb Z$
which is defined naturally by induction on $n$.
\qed\endpf

Compare the following Proposition 3 with
Weil \cite{17, Ch. II \S5 Theorem 3}.

\th Proposition 3

Let $K$ be an  $n$-dimensional local field,
and let $V$ be a vector space over $K$ of finite dimension,
$V'=\Hom_K(V,K)$. Then
\Roster :
\Item{(1)} The abelian group object $\U V$ of $\Cal F_n$ which represents
the functor
$X\to V\otimes_K [X,\U K]$ belongs to $\Cal A_n$.

\Item{(2)} $[\U K,\Bbb Q/\Bbb Z]_{\Cal A_\infty}$ is one-dimensional
with respect to the natural $K$-module structure
and its non-zero element induces due to Lemma 3 (2) an isomorphism
$\U{V'}\simeq D(\U V)$.
\endRoster 
\endth

\HH 3. Properties of the ring of $\U K$-valued morphisms

\HHH 3.1. Multiplicative groups of certain complete rings

\th Proposition 4

Let $A$ be a ring and let $\pi$ be a non-zero element of $A$ such that
$A=\varprojlim A/\pi^nA$.
Let $R=A/\pi A$ and $B=A[\pi^{-1}]$.
Assume that at least one of the following two conditions is satisfied.
\Roster 
\Item{(i)} $R$ is reduced (i.e. having no nilpotent elements except zero)
and there is a ring homomorphism $s\colon R\to A$ such that the composite
$R@>s>> A @>>> A/\pi A$ is the identity.

\Item{(ii)} For a prime $p$ the ring $R$ is annihilated by $p$
and via the homomorphism $R\to R$, $x\mapsto x^p$ the latter $R$ is a finitely generated
projective module over the former $R$.
\endRoster

Then we have
$$B^*\simeq A^*\times \Gamma(\Spec(R),\Bbb Z)$$
where $\Gamma(\Spec(R),\Bbb Z)$ is the group of global sections of the constant sheaf
$\Bbb Z$ on $\Spec(R)$ with Zariski topology.
The isomorphism is given by the homomorphism
of sheaves $\Bbb Z \to \Cal O_{\Spec(R)}^*$, $1\mapsto \pi$,
the map
$$\Gamma(\Spec(R),\Bbb Z)\simeq \Gamma(\Spec(A),\Bbb Z) 
\to \Gamma(\Spec(B),\Bbb Z)$$
and the inclusion map $A^*\to B^*$.
\endth
\pf Proof

Let $\text{\rm Aff}_R$ be the category of affine schemes
over $R$.
In  case (i) let $\Cal C=\text{\rm Aff}_R$.
In case (ii) let $\Cal C$ be the category
of all affine schemes $\Spec(R')$ over $R$ such that the map
$$R^{(p)}\otimes_R R'\to R',\quad x\otimes y\mapsto xy^p$$
(cf.\ the proof of Proposition 2) is bijective.
Then in case (ii) every finite inverse limit and finite sum
exists in $\Cal C$
and coincides with that taken in $\text{\rm Aff}_R$.
Furthermore, in this case the inclusion functor
$\Cal C\to \text{\rm Aff}_R$ has a right adjoint.
Indeed, for any affine scheme $X$ over $R$ the corresponding object
in $\Cal C$ is
$\varprojlim X_i$ where $X_i$ is the Weil restriction of $X$ 
with respect to the homomorphism $R\to R$, $x\mapsto x^{p^i}$.

Let $\U R$ be the ring object of $\Cal C$ which represents
the functor
$X\to \Gamma(X,\Cal O_X)$,
and let $\U{R^*}$ be the object which represents the functor
$X\to [X,\U R]^*$,
and $\U 0$ be the final object $e$ regarded as a closed subscheme of $\U R$
via the zero morphism $e\to \U R$.

\th Lemma 4

Let $X$ be an object of $\Cal C$ and assume that $X$ is reduced as a scheme
(this condition is always satisfied in case (ii)).
Let $\theta\colon X\to \U R$ be a morphism of $\Cal C$.
If $\theta^{-1}(\U{R^*})$ is a closed subscheme of $X$, 
then $X$ is the direct sum of
$\theta^{-1}(\U{R^*})$ and $\theta^{-1}(\U 0)$
{\rm(}where the inverse image notation are used for the fibre product{\rm)}.
\endth

The group $B^*$ is generated by elements $x$ of $A$ such that $\pi^n\in Ax$ for some
$n\ge 0$.
In case (i) let $\U{A/\pi^{n+1}A}$ be the ring object of $\Cal C$
which represents the functor
$X\to A/\pi^{n+1}A\otimes_R [X,\U R]$
where
$A/\pi^{n+1}A$ is viewed as an $R$-ring via a fixed section $s$.
In case (ii) we get a ring object 
$\U{A/\pi^{n+1}A}$ of $\Cal C$ by Proposition 2 (4).

In both cases there are morhisms $\theta_i\colon \U R\to \U{A/\pi^{n+1}A}$
($0\le i\le n$) in $\Cal C$ 
such that the morphism 
$$\U R\times\dots\times\U R\to \U{A/\pi^{n+1}A},\quad (x_0,\dots,x_n)\mapsto \sum_{i=0}^n \theta_i(x_i)\pi^i
$$
is an isomorphism.

Now assume $xy=\pi^n$ for some $x,y\in A$ and take elements
$x_i,y_i\in R=[e,\U R]$ ($0\le i\le n$) such that
$$x\mod \pi^{n+1} = \sum_{i=0}^n \theta_i(x_i)\pi^i,\quad 
y\mod \pi^{n+1} = \sum_{i=0}^n \theta_i(y_i)\pi^i.
$$
An easy computation shows that
for every $r=0,\dots, n$
$$\bigl(\bigcap_{i=0}^{r-1}x_i^{-1}(\U 0)\bigr)\bigcap x_r^{-1}(\U{R^*})=
\bigl(\bigcap_{i=0}^{r-1}x_i^{-1}(\U 0)\bigr)\bigcap 
\bigl(\bigcap_{i=0}^{n-r-1}y_i^{-1}(\U 0)\bigr).$$
By Lemma 4 and induction on $r$ we deduce that 
$e=\Spec(R)$ is the direct sum of the closed open
subschemes $\bigl(\cap_{i=0}^{r-1}x_i^{-1}(\U 0)\bigr)\cap x_r^{-1}(\U{R^*})$
on which the restriction of $x$ has the form $a\pi^r$ for an invertible
element $a\in A$.
\qed\endpf

\HHH 3.2. Properties of the ring $[X,\U K]$

\phantom{}\smallskip\par

Results of this subsection will be used in section 4.

\df Definition 4

For an object $X$ of $\Cal F_\infty$ and a set $S$ let
$$\lcf(X,S)=\inlim_{I} \, [X,I]$$
where $I$ runs over all finite subsets of $S$
(considering each $I$ as an object of $\Cal F_0\subset \Cal F_\infty$).
\enddf

\th Lemma 5

Let $K$ be an $n$-dimensional local field
and let $X$ be an object of $\Cal F_\infty$.
Then:
\Roster 
\Item{(1)} The ring $[X,\U K]$ is reduced.

\Item{(2)} For every set $S$ there is a canonical bijection
$$\lcf(X,S)\iss \Gamma(\Spec([X,\U K]),S)$$
where $S$ on the right hand side is regarded as a constant
sheaf on $\Spec([X,\U K])$.
\endRoster 
\endth
\pf Proof of (2)

If $I$ is a finite set and $\theta\colon X\to I$ is a morphism of
$\Cal F_\infty$ then $X$ is the direct sum of the objects
$\theta^{-1}(i)=X\times_I \{i\}$ in $\Cal F_\infty$ ($i\in I$).
Hence we get the canonical map of (2).
To prove its bijectivity we may assume $S=\{0,1\}$.
Note that $\Gamma(\Spec(R),\{0,1\})$ is the set of idempotents
in $R$ for any ring $R$.
We may assume that $X$ is an object of $\pro(\Cal F_{n-1})$.

Let $k_{n-1}$ be the residue field of $k_n=K$.
Then $$\Gamma(\Spec([X,\U K]),\{0,1\})\simeq
\Gamma(\Spec([X,\U{k_{n-1}}]),\{0,1\})$$
by (1) applied to the ring $[X,\U{k_{n-1}}]$.
\qed\endpf

\th Lemma 6

Let $K$ be an $n$-dimensional local field of characteristic $p>0$.
Let 

\noindent $k_0,\dots,k_n$ be as in the introduction.
For each $i=1,\dots, n$ let $\pi_i$ be a lifting to $K$ of
a prime element of $k_i$.
Then for each object $X$ of $\Cal F_\infty$
$[X,\U K]^*$ is generated by the subgroups
$$[X,\U{K^p(\pi^{(s)})}]^*$$ 
where $s$ runs over
all functions $\{1,\dots,n\}\to \{0,1,\dots,p-1\}$
and $\pi^{(s)}$ denotes $\pi_1^{s(1)}\dots \pi_n^{s(n)}$,
$\U{K^p(\pi^{(s)})}$ is the subring object of $\U K$ corresponding to
$K^p(\pi^{(s)})$, i.e.
$$[X,\U{K^p(\pi^{(s)})}]=K^p(\pi^{(s)})\otimes_{K^p}[X,\U K].$$
\endth
\pf Proof

Indeed,  Proposition 4
and induction on $n$ yield 
morphisms $$\theta^{(s)}\colon \U{K^*}\to \U{K^p(\pi^{(s)})^*}$$ 
such that the product of all $\theta^{(s)}$ in $\U{K^*}$
is the identity morphism $\U{K^*}\to \U{K^*}$.
\qed\endpf

The following similar result is also proved by induction on $n$.

\th Lemma 7

Let $K,k_0$ and $(\pi_i)_{1\le i\le n}$ be as in Lemma 6.
Then there exists a morphism of $\Cal A_\infty$ 

\noindent (cf.\ section 2)
$$(\theta_1,\theta_2)\colon \U{\Omega_K^n}\to 
\U{\Omega_K^n}\times \U{k_0}$$
such that
$$x=(1-\text{\tenrm C})\theta_1(x)+\theta_2(x)d\pi_1/\pi_1\wedge\dots \wedge d\pi_n/\pi_n$$
for every object $X$ of $\Cal F_\infty$ and for every $x\in [X,\U{\Omega_K^n}]$
where $\U{\Omega_K^n}$ is the object which represents the functor
$X\to \Omega_K^n\otimes_K[X,\U K]$ and $\text{\tenrm C}$ denotes
the Cartier operator {\rm(}\cite{4}, or  see 4.2 in Part I
for the definition{\rm)}.
\endth

Generalize the Milnor $K$-groups as follows.

\df Definition 5

For a ring $R$ let $\Gamma_0(R)=\Gamma(\Spec(R),\Bbb Z)$.
The morphism of sheaves
$$\Bbb Z\times \Cal O_{\Spec(R)}^*\to \Cal O_{\Spec(R)}^*,
\quad (n,x)\mapsto x^n$$
determines the $\Gamma_0(R)$-module structure on $R^*$.
Put $\Gamma_1(R)=R^*$ and for $q\ge 2$ put
$$\Gamma_q(R)=\otimes_{\Gamma_0(R)}^q \Gamma_1(R)/J_q$$
where 
$\otimes_{\Gamma_0(R)}^q \Gamma_1(R)$ is the $q$th tensor power of
$\Gamma_1(R)$ over $\Gamma_0(R)$ and
$J_q$ is the subgroup of the tensor power generated by elements
$x_1\otimes\dots\otimes x_q$ which satisfy $x_i+x_j=1$ or $x_i+x_j=0$ for some $i\not=j$.
An element $x_1\otimes\dots\otimes x_q\mod J_q$ will be denoted by $\{x_1,\dots, x_q\}$.
\enddf

Note that $\Gamma_q(k)=K_q(k)$ for each field $k$
and $\Gamma_q(R_1\times R_2)\simeq \Gamma_q(R_1)\times
\Gamma_q(R_2)$ for rings $R_1$, $R_2$.

\th Lemma 8

In one of the following two cases
\Roster 
\Item{(i)} $A,R,B,\pi$ as in Proposition 4 

\Item{(ii)} an $n$-dimensional local field $K$,
an object $X$ of $\Cal F_\infty$,
$A=[X,\U{\Cal O_K}]$, 

\Item{} $R=[X,\U{k_{n-1}}]$, $B=[X,\U{K}]$,
\endRoster 

\noindent let $U_i\Gamma_q(B)$ be the subgroup of $\Gamma_q(B)$ generated
by elements $\{1+\pi^ix,y_1,\dots,y_{q-1}\}$ such that
$x\in A$, $y_j\in B^*$, $q,i\ge 1$.

Then: 
\Roster 
\Item{(1)} There is a homomorphism
$\rho_0^q\colon \Gamma_q(R)\to \Gamma_q(B)/U_1\Gamma_q(B)$
such that
$$\rho_0^q(\{x_1,\dots,x_q\})=\{\wt{x_1},\dots,\wt{x_q}\}\mod U_1\Gamma_q(B)$$
where $\wt{x_i}\in A$ is a representative of$x_i$.
In case (i) {\rm(}resp.\ (ii){\rm)} the induced map
$$\Gamma_q(R)+\Gamma_{q-1}(R)\to \Gamma_q(B)/U_1\Gamma_q(B),
\quad (x,y)\mapsto \rho_0^q(x)+\{\rho_0^{q-1}(y),\pi\}
$$
{\rm(}resp.\
$$\aligned
&\Gamma_q(R)/m+\Gamma_{q-1}(R)/m\to \Gamma_q(B)/(U_1\Gamma_q(B)+m\Gamma_q(B)),\\
& (x,y)\mapsto \rho_0^q(x)+\{\rho_0^{q-1}(y),\pi\}\text{\rm)}
\endaligned 
$$
is bijective {\rm(}resp.\ bijective for every non-zero integer $m${\rm)}.
\it 
\Item{(2)} If $m$ is an integer invertible in $R$ then $U_1\Gamma_q(B)$ is $m$-divisible.

\Item{(3)} In case (i) assume that $R$ is additively generated by $R^*$.
In case (ii) assume that $\chr(k_{n-1})=p>0$.
Then there exists a unique homomorphism
$$\rho_i^q\colon \Omega_R^{q-1}\to U_i\Gamma_q(B)/U_{i+1}\Gamma_q(B)$$ such that
$$\rho_i^q(x dy_1/y_1\wedge\dots\wedge dy_{q-1}/y_{q-1})=
\{1+\wt{x}\pi^i, \wt{y_1},\dots,\wt{y_{q-1}}\} \mod U_{i+1}\Gamma_q(B)$$
for every $x\in R$, $y_1,\dots, y_{q-1}\in R^*$.
The induced map
$$\Omega_R^{q-1}\oplus \Omega_R^{q-2}\to U_i\Gamma_q(B)/U_{i+1}\Gamma_q(B),
\quad (x,y)\mapsto \rho_i^q(x)+\{\rho_i^{q-1}(y),\pi\}$$
is surjective.
If $i$ is invertible in $R$ then the homomorphism $\rho_i^q$ is surjective.
\endRoster
\endth
\pf Proof

In case (i) these results follow from Proposition 4 by Bass--Tate's method
\cite{2, Proposition 4.3} for (1),
Bloch's method \cite{3, \S3} for (3)
and by writing down the kernel of
$R\otimes R^*\to \Omega_R^1$, $x\otimes y\mapsto xdy/y$ as in \cite{9, \S1 Lemma 5}.

If $X$ is an object of $\pro(\Cal F_{n-1})$ then case (ii)
is a special case of (i) except 
$n=1$ and $k_{0}=\Bbb F_2$ where $[X,\U{k_0}]$ is not generated by
$[X,\U{k_0}]^*$ in general.
But in this exceptional case it is easy to check directly all the assertions.

For an arbitrary $X$ we present here only the proof of (3) because
the proof of (1) is rather similar.

Put $k=k_{n-1}$.
For the existence of $\rho_i^q$ it suffices to consider the cases where 
$X=\U{\Omega_{k}^{q-1}}$ and
$X=\U k\times\prod^{q-1} \U{k^*}$ ($\prod^r Y$ denotes the product of
$r$ copies of $Y$).
Note that these objects are in $\pro(\Cal F_{n-1})$ since
 $[X,\U{\Omega_{k}^q}]=\Omega_{[X,\U k]}^q$ for any $X$ and $q$.

The uniqueness follows from the fact that
$[X,\U{\Omega_{k}^{q-1}}]$ is generated by elements of the form 
$x d c_1/c_1\wedge\dots\wedge d c_{q-1}/c_{q-1}$
such that
$x\in [X,\U{k}]$ and $c_1,\dots,c_{q-1}\in k^*$.

To prove the surjectivity we may assume $X=(1+\pi^i\U{\Cal O_K})\times
\prod^{q-1}\U{K^*}$
and it suffices to prove in this case that the typical element
in $U_i\Gamma_q(B)/U_{i+1}\Gamma_q(B)$ belongs to
the image of the homomorphism introduced in (3).
Let $\U{U_K}$ be the object of $\Cal F_n$ which represents
the functor
$X\to [X,\U{\Cal O_K}]^*$.
By Proposition 4 there exist 

\noindent morphisms
$\theta_1\colon \U{K^*}\to \coprod_{i=0}^{p-1} \U{U_K}\pi^i$
(the direct sum in $\Cal F_n$) 
and
$\theta_2\colon \U{K^*}\to \U{K^*}$

\noindent such that $x=\theta_1(x)\theta_2(x)^p$ for each $X$ in $\Cal F_\infty$ and 
each $x\in [X,\U{K^*}]$
(in the proof of (1) $p$ is replaced by $m$).
Since $\coprod_{i=0}^{p-1} \U{U_K}\pi^i$ belongs to $\pro(\Cal F_{n-1})$
and 

\noindent $(1+\pi^i[X,\U{\Cal O_K}])^p\subset 1+\pi^{i+1}[X,\U{\Cal O_K}]$ we are reduced to the case where $X$ is an object of $\pro(\Cal F_{n-1})$.
\qed\endpf

\HH 4. Norm groups

\phantom{}
\par

In this section we prove Theorem 3 and Proposition 1.
In subsection 4.1 we reduce these results to Proposition 6.

\HHH 4.1. Reduction steps

\df Definition 6

Let $k$ be a field
and let $H\colon \Cal E(k)\to\text{\rm Ab}$ be a functor
such that 

\noindent $\inlim_{k'\in \Cal E(k)}\, H(k')=0$.
Let $w\in H(k)$ (cf.\ Introduction).
For a ring $R$ over $k$ and $q\ge 1$ define the subgroup
$N_q(w,R)$ (resp.\ $L_q(w,R)$) of $\Gamma_q(R)$ as follows.

An element $x$ belongs to  $N_q(w,R)$ (resp.\ $L_q(w,R)$)
if and only if there exist 
\Roster 
\Item{} a finite set $J$ and element $0\in J$,

\Item{} a map $f\colon J\to J$ such that for some $n\ge 0$
the $n$th iteration $f^n$ with respect to the composite
is a constant map with value $0$,

\Item{} and a family $(E_j,x_j)_{j\in J}$
($E_j\in \Cal E(k)$), $x_j\in \Gamma_q(E_j\otimes_k R)$)
satisfying the following conditions: 
\ItemItem{(i)} $E_0=k$ and $x_0=x$.
\ItemItem{(ii)} $E_{f(j)}\subset E_j$ for every $j\in J$.
\ItemItem{(iii)} Let $j\in f(J)$. Then there exists a family
$(y_t,z_t)_{t\in f^{-1}(j)}$  
\ItemItem{} 
($y_t\in (E_t\otimes_k R)^*$, $z_t\in\Gamma_{q-1}(
E_j\otimes_k R)$) such that
$x_t=\{y_t, z_t\}$ for all $t\in f^{-1}(j)$ and
$$x_j=\sum_{t\in f^{-1}(j)}\{N_{E_t\otimes_k R/E_j\otimes_kR}(y_t),z_t\}$$
where
$N_{E_t\otimes_k R/E_j\otimes_kR}$ denotes the norm homomorphism
$$(E_t\otimes_k R)^*\to (E_j\otimes_kR)^*.$$
\ItemItem{(iv)} If $j\in J\setminus f(J)$ then
$w$ belongs to the kernel of
$H(k)\to H(E_j)$
\ItemItem{}(resp.\ then one of the following two assertions is valid:
\ItemItem{\qquad(a)} $w$ belongs to the kernel of
$H(k)\to H(E_j)$,
\ItemItem{\qquad(b)} $x_j$ belongs to the image of
$\Gamma(\Spec(E_j\otimes_kR), K_q(E_j))\to \Gamma_q(E_j\otimes_k R)$,
where $K_q(E_j)$ 
denotes the constant sheaf on
$\Spec(E_j\otimes_k R)$ defined by the set $K_q(E_j)$).
\endRoster
\enddf

\rk Remark

If the groups $\Gamma_q(E_j\otimes_k R)$ have a suitable ``norm'' homomorphism
then $x$ is the sum of the ``norms'' of $x_j$ such that 
$f^{-1}(j)=\emptyset$.
In particular, in  the case where $R=k$ we get
$N_q(w,k)\subset N_q(w)$ and $N_1(w,k)=N_1(w)$.
\endrk

\df Definition 7

For a field $k$ let
$[\Cal E(k),\text{\rm Ab}]$ be the abelian category of all functors
$\Cal E(k)\to \text{\rm Ab}$.
\Roster 
\Item{(1)} For $q\ge 0$ let $\Cal N_{q,k}$ denote the full subcategory
of $[\Cal E(k),\text{\rm Ab}]$
consisting of functors $H$ such that
$\inlim_{k'\in \Cal E(k)}\, H(k')=0$
and such that for every $k'\in \Cal E(k)$, 
$w\in H(k')$ the norm group
$N_q(w)$ coincides with $K_q(k')$.
Here $N_q(w)$ is defined with respect to the functor $\Cal E(k')\to \text{\rm Ab}$.

\Item{(2)} If $K$ is an $n$-dimensional local field and $q\ge 1$, 
let $\U{\Cal N}_{q,K}$ (resp.\ $\U{\Cal L}_{q,K}$)
denote the full subcategory of
$[\Cal E(K),\text{\rm Ab}]$
consisting of functors $H$ such that
$$\inlim_{K'\in \Cal E(K)} \, H(K')=0$$
and such that for every $K'\in \Cal E(K)$, 
$w\in H(K')$ and every object $X$ of $\Cal F_\infty$
the group $N_q(w,[X,\U{K'}])$
(resp.\ $L_q(w,[X,\U{K'}])$)
coincides with $\Gamma_q([X,\U{K'}])$.
\endRoster
\enddf

\th Lemma 9

Let $K$ be an $n$-dimensional local field and let $H$ be an object of
$\U{\Cal L}_{q,K}$.
Then for every $w\in H(K)$
the group $N_q(w)$ is an open subgroup of $K_q(K)$ of finite index.
\endth

\pf Proof

Consider the case where 
$X=\prod^q \U{K^*}$.
We can take a system
$(E_j,x_j)_{j\in J}$
as in Definition 6 such that
$E_0=K$, $x_0$ is the canonical element
in $\Gamma_q([X,\U K])$
and such that
if $j\not\in f(J)$ and $w\not\in \kr(H(K)\to H(E_j))$
then $x_j$ is the image of an element $\theta_j$ of
$\lcf(X,K_q(E_j))$.
Let $\theta\in \lcf(X,K_q(K)/N_q(w))$ be the sum of
$N_{E_j/K}\circ \theta_j\mod N_q(w)$.
Then the canonical map
$[e,X]=\prod^q K^*\to K_q(K)/N_q(w)$ comes from $\theta$.
\qed\endpf

\df Definition 8

Let $k$ be a field.
A collection $\{\Cal C_{k'}\}_{k'\in \Cal E(k)}$ of full
subcategories $\Cal C_{k'}$ of
$[\Cal E(k'),\text{\rm Ab}]$
is called {\it admissible}  if and only if
it satisfies conditions (i) -- (iii) below.
\Roster 
\Item{(i)} Let $E\in \Cal E(k)$.
Then every subobject, quotient object, extension and filtered inductive limit
(in the category of $[\Cal E(E),\text{\rm Ab}]$) of objects of $\Cal C_E$
belongs to $\Cal C_E$.

\Item{(ii)} Let $E,E'\in\Cal E(k)$ and
$E\subset E'$.
If $H$ is in $\Cal C_E$ then the composite functor
$\Cal E(E')@>>> \Cal E(E)@>H>> \text{\rm Ab}$ is in $\Cal C_{E'}$.

\Item{(iii)} Let $E\in \Cal E(k)$ and $H$ is in $[\Cal E(E),\text{\rm Ab}]$.
Then $H$ is in $\Cal C_E$ if  conditions (a) and (b) below are satisfied for a prime $p$.

\ItemItem{(a)} For some $E'\in \Cal E(E)$ such that $|E':E|$ is prime to $p$
the composite functor
$(E')@>>> (E)@>H>> \text{\rm Ab}$ is in $\Cal C_{E'}$.

\ItemItem{(b)} Let $q$ be a prime number distinct from $p$
and let $S$ be a direct subordered set of $\Cal E(E)$.
If the degree of every finite extension of the field
$\inlim_{E'\in S} \, E'$ is a power of $p$
then $\inlim_{E'\in S}\, H(E')=0$. 
\endRoster 
\enddf

\th Lemma 10

\Roster 
\Item{(1)} For each field $k$ and $q$ the collection
$\{\Cal N_{q,k'}\}_{k'\in \Cal E(k)}$
is admissible.
If $K$ is an $n$-dimensional local field then the collections
$\{\U{\Cal N}_{q,k'}\}_{k'\in \Cal E(k)}$ and
$\{\U{\Cal L}_{q,k'}\}_{k'\in \Cal E(k)}$
are admissible.

\Item{(2)} Let $k$ be a field. Assume that a collection
$\{\Cal C_{k'}\}_{k'\in \Cal E(k)}$ is admissible.
Let $r\ge 1$ and for every prime $p$ there exist
$E\in \Cal E(k)$ such that
$|E:k|$ is prime to $p$ and such that the functor
$H^r(\,\,\,,\Bbb Z/p^r)\colon \Cal E(E)\to\text{\rm Ab}$ is in $\Cal C_E$.
Then for each $k'\in \Cal E(k)$, each discrete torsion
abelian group $M$ endowed with a continuous action of
$\Gal({k'}^{\sep}/k')$ and each $i\ge r$
the functor $$H^i(\,\,\,,M)\colon \Cal E(k')\to\text{\rm Ab}$$
is in $\Cal C_{k'}$.
\endRoster
\endth

\df Definition 9

For a field $k$, $r\ge 0$ and a non-zero integer $m$ define the group $H_m^r(k)$ as follows.

If $\chr(k)=0$ let
$$H^r_m(k)=H^r(k, \mu_m^{\otimes(r-1)}).$$
If $\chr(k)=p>0$ and $m=m'p^i$ where $m'$ is prime to $p$ and $i\ge 0$ let
$$H_{m}^r(k)=H_{m'}^r(k,\mu_{m'}^{\otimes(r-1)})\oplus\text{\rm coker}
(F-1\colon C_i^{r-1}(k)\to C_i^{r-1}(k)/\{C_i^{r-2}(k),T\})$$
(where 
$C_i^{\cdot}$ is the group defined in \cite{3, Ch.II,\S 7}, 
$C_i^r=0$ for $r<0$).
\enddf

{\Em {By the above results it suffices for the proof of Theorem 3
to prove the following Proposition 5 in  the case where $m$ is a prime number.}}

\th Proposition 5

Let $K$ be an $n$-dimensional local field.
Let $q,r\ge 1$ and let $m$ be a non-zero integer.
Then the functor $H_m^r\colon \Cal E(K)\to\text{\rm Ab}$
is in $\U{\Cal L}_{q,K}$ if $q+r=n+1$ and
in 
$\U{\Cal N}_{q,K}$ if $q+r> n+1$.
\endth

\smallskip

{\Em{Now we begin the proofs of Proposition 1 and Proposition 5.}}

\df Definition 10

Let $K$ be a complete discrete valuation field,
$r\ge 0$ and $m$ be a non-zero integer.
\Roster

\Item{(1)} Let
$H^r_{m,\ur}$ and $H^r_m/H^r_{m,\ur}$ be the functors
$\Cal E(K)\to \text{\rm Ab}$:
$$\aligned
& H^r_{m,\ur}(K')=\kr(H^r_{m}(K')\to H^r_{m}(K'_{\ur})), \\
&(H^r_m/H^r_{m,\ur})(K')=H_m^r(K')/H^r_{m,\ur}(K')
\endaligned
$$
where $K'_{\ur}$ is the maximal unramified extension of $K'$.

\Item{(2)} Let $I_m^r$ (resp.\ $J_m^r$)
be the functor
$\Cal E(K)\to \text{\rm Ab}$ such that
$I_m^r(K')=H_m^r(k')$ (resp.\  $J_m^r(K')=H_m^r(k')$)
where $k'$ is the residue field of $K'$
and such  that
the homomorphism
$I_m^r(K')\to I_m^r(K'')$ 
(resp.\ $J_m^r(K')\to J_m^r(K'')$)
for $K'\subset K''$ is
$j_{k''/k'}$ (resp.\ $e(K''|K')j_{k''/k'}$)
where $k''$ is the residue field of $K''$,
$j_{k''/k'}$ is the canonical homomorphism induced by the inclusion
$k'\subset k''$ and $e(K''|K')$ is the index of ramification
of $K''/K'$.
\endRoster 
\enddf

\th Lemma 11

Let $K$ and $m$ be as in Definition 10.
\Roster 
\Item{(1)} For $r\ge 1$ there exists an exact sequence
of functors
$$0\to I_m^r\to H_{m,\ur}^r\to J_m^{r-1}\to 0.$$

\Item{(2)} $J_m^r$ is in $\Cal N_{1,K}$ for every $r\ge 0$.

\Item{(3)} Let $q,r\ge 1$. Then $I_m^r$ is in $\Cal N_{q,K}$
if and only if $H_m^r\colon \Cal E(k)\to \text{\rm Ab}$
is in $\Cal N_{q-1,k}$ where $k$ is the residue field of $K$.
\endRoster 
\endth
\pf Proof

The assertion (1) follows from \cite{11}.
The assertion (3) follows from the facts that
$1+\Cal M_K\subset N_{L/K}(L^*)$ for every unramified
extension $L$ of $K$
and that
there exists a canonical split exact sequence
$$0\to K_q(k)\to K_q(K)/U_1K_q(K)\to K_{q-1}(k)\to 0.\eqno{\qed}$$
\endpf

The following proposition will be proved in 4.4.

\th Proposition 6

Let $K$ be a complete discrete valuation field
with residue field $k$.
Let $q,r\ge 1$ and $m$ be a non-zero integer.
Assume that $[k:k^p|\le p^{q+r-2}$ 
if $\chr(k)=p>0$. Then:
\Roster 
\Item{(1)} $H_m^r/H_{m,\ur}^r$ is in $\Cal N_{q,K}$.

\Item{(2)} If $K$ is an $n$-dimensional local field
with $n\ge 1$ then
 $H_m^r/H_{m,\ur}^r$ is in $\U{\Cal N}_{q,K}$.
\endRoster 
\endth

{\Em {Proposition 1 follows from this proposition
by Lemma 10 and Lemma 11}} 
(note that if $\chr(k)=p>0$ and $i\ge 0$ then
$H_{p^i}^r(k)$ is isomorphic to
$\kr(p^i\colon H^r(k)\to H^r(k))$ as it follows from \cite{11}).

\th Lemma 12

Let $K$ be an $n$-dimensional local field
and let $X$ be an object of $\Cal F_\infty$.
Consider the following cases.
\Roster 
\Item{(i)} $q>n+1$ and $m$ is a non-zero integer.

\Item{(ii)} $q=n+1$, $\chr(K)=p>0$ and $m$ is a power of $p$.

\Item{(iii)} $q=n+1$ and $m$ is a non-zero integer.
\medskip 
Let $x\in \Gamma_q([X,\U K]$.
Then in cases (i) and (ii)
{\rm(}resp.\ in case (iii){\rm)} there exist a triple
$(J,0,f)$ and a family
$(E_j,x_j)_{j\in J}$ which satisfy
all the conditions in Definition 6 with $k=K$
except condition (iv),
and which satisfy the following condition:

\Item{(iv)'} If $j\in J\setminus f(J)$ then
$x_j\in m\Gamma_q([X,\U{E_j}])$
\Item{}{\rm(}resp.\ $x_j$ belongs to $m\Gamma_q([X,\U{E_j}])$ 
\Item{} or
to the image of $\lcf(X,K_q(E_j))\to \Gamma_q([X,\U{E_j}])${\rm)}.
\endRoster 
\endth

\th Corollary

Let $K$ be an $n$-dimensional local field.
Then  $mK_{n+1}(K)$ is an open subgroup of finite index of  $K_{n+1}(K)$ for  every non-zero integer $m$.
\endth

This corollary follows from  case (iii) above by the argument in the proof
of Lemma 9.

\smallskip

\pf Proof of Lemma 12

We may assume that $m$ is a prime number.

First we consider case (ii).
By Lemma 6 we may assume that there are elements
$b_1,\dots,b_{n+1}\in [X,\U K]^*$ and
$c_1,\dots, c_{n+1}\in K^*$ such that
$x=\{b_1,\dots,b_{n+1}\}$ and
$b_i\in [X,\U{K^p(c_i)}]^*$ for each $i$.
We may assume that 
$|K^p(c_1,\dots,c_r):K^p|=p^r$ and $c_{r+1}\in K^p(c_1,\dots,c_r)$
for some $r\le n$.
Let $J=\{0,1,\dots, r\}$,
and define $f\colon J\to J$ by
$f(j)=j-1$ for $j\ge 1$ and $f(0)=0$.
Put $E_j=K(c_1^{1/p},\dots,c_j^{1/p})$
and $x_j=\{b_1^{1/p},\dots,b_j^{1/p},b_{j+1},\dots,b_{n+1}\}$.
Then

\noindent $x_r=p\{b_1^{1/p},\dots, b_{r+1}^{1/p},b_{r+2},\dots,b_{n+1}\}$
in $\Gamma_{n+1}([X,\U{E_r}])$.

Next we consider cases (i) and (iii).
If $K$ is a finite field
then the assertion for (i) follows from Lemma 13 below
and the assertion for (iii) is trivial.
Assume $n\ge 1$ and let $k$ be the residue field of $K$.
By induction on $n$
Lemma 8 (1) (2) and case (ii) of Lemma 12 show that
we may assume
$x\in U_1\Gamma_q([X,\U K])$, $\chr(K)=0$ and $m=\chr(k)=p>0$.
Furthermore
we may assume that $K$ contains a primitive $p$th root $\zeta$ of 1.
Let $e_K=v_K(p)$ and let $\pi$ be a prime
element of $K$.
Then 
$$
U_i\Gamma_q([X,\U{\Cal O_K}])\subset pU_1\Gamma_q([X,\U{\Cal O_K}]),
\quad\text{\rm  
if $i>pe_K/(p-1)$.}
$$
From this and Lemma 8 (3)
(and a computation of the map $x\mapsto x^p$ on $U_1\Gamma_q([X,\U{\Cal O_K}])$)
it follows that
$U_1\Gamma_q([X,\U{K}])$ is $p$-divisible if $q>n+1$
and that there is a surjective homomorphism
$$
\aligned
&[X,\U{\Omega_k^{n-1}}]/(1-\text{\tenrm C})[X,\U{\Omega_k^{n-1}}]\to
U_1\Gamma_{n+1}([X,\U{K}])/pU_1\Gamma_{n+1}([X,\U{K}]),
\\
&xdy_1/y_1\wedge\dots\wedge dy_{n-1}/y_{n-1}\mapsto
\{1+\wt{x}(\zeta-1)^p,\wt{y_1},\dots,\wt{y_{n_1}},\pi\}
\endaligned
$$
where $\text{\tenrm C}$ is the Cartier operator.
By Lemma 7
$$[X,\U{\Omega_k^{n-1}}]/(1-\text{\tenrm C})[X,\U{\Omega_k^{n-1}}]
=\lcf(X,\Omega_k^{n-1}/(1-\text{\tenrm C})\Omega_k^{n-1}).\eqno{\qed}$$
\endpf

\th Lemma 13

Let $K$ be a finite field and let $X$ be an object of $\Cal F_\infty$.
Then
\Roster 
\Item{(1)} $\Gamma_q[X,\U K]=0$ for $q\ge 2$.

\Item{(2)} For every finite extension $L$ of $K$
the norm homomorphism
$[X,\U L]^*\to [X,\U K]^*$ is surjective.
\endRoster \endth
\pf Proof

Follows from Lemma 5 (2).
\qed\endpf

\pf Proof of Proposition 5 assuming Proposition 6

If $K$ is a finite field,  the assertion of Proposition 5 follows from Lemma 13.

Let $n\ge 1$. Let $k$ be the residue field of $K$.
Let $I_m^r$ and $J_m^r$ be as in Definition 10.
Assume $q+r=n+1$ 
(resp.\ $q+r>n+1$).
Using Lemma 8 (1) and the fact that
$$U_1\Gamma_q([X,\U K])\subset N_{L/K}\Gamma_q([X,\U L])$$
for every unramified extension $L/K$ we can deduce that
$I_m^r$ is in $\U{\Cal L}_{q,K}$
(resp.\ $\U{\Cal N}_{q,K}$)
from the induction hypothesis 
$H_m^r\colon \Cal E(k)\to \text{\rm Ab}$ is in
$\U{\Cal L}_{q-1,k}$
(resp.\ $\U{\Cal N}_{q-1,k}$).
We can deduce $J_m^{r-1}$ is in $\U{\Cal L}_{q,K}$
(resp.\ $\U{\Cal N}_{q,K}$) from the hypothesis
$H^{r-1}_m\colon \Cal E(k)\to \text{\rm Ab}$ is in 
 $\U{\Cal L}_{q,k}$
(resp.\ $\U{\Cal N}_{q,k}$).
Thus $H_{m,\ur}^r$ is in  $\U{\Cal L}_{q,K}$
(resp.\ $\U{\Cal N}_{q,K}$).
\qed\endpf

\HHH 4.2. Proof of Proposition 6

\phantom{}\smallskip\par

Let $k$ be a field and let $m$ be a non-zero integer.
Then
$\oplus_{r\ge 0}H_m^r(k)$ (cf.\ Definition 9)
has a natural right
$\oplus_{q\ge0} K_q(k)$-module structure
(if $m$ is invertible in $k$ this structure is defined by the cohomological symbol
$h_{m,k}^q\colon K_q(k)/m\to H^q(k,\mu_m^{\otimes q})$
and the cup-product, cf.\ \cite{9, \S3.1}).
We denote the product in this structure
by
$\{w,a\}$ 

\noindent ($a\in \oplus_{q\ge0} K_q(k)$, $w\in \oplus_{r\ge 0}H_m^r(k)$).

\df Definition 11

Let $K$ be a complete discrete valuation field
with residue field $k$ such that $\chr(k)=p>0$.
Let $r\ge 1$.
We call an element $w$ of $H_p^r(K)$
{\it standard} if and only if
$w$ is in  one of the following forms (i) or (ii).
\Roster 
\Item{(i)} 
$w=\{\chi,a_1,\dots,a_{r-1}\}$ where $\chi$ is an element
of $H_p^1(K)$ corresponding to a totally ramified cyclic extension
of $K$ of degree $p$,
and $a_1,\dots,a_{r-1}$ are elements of $\Cal O_K^*$ such that
$$|k^p(\Overline{a_1},\dots,\Overline{a_{r-1}}):k^p|=p^{r-1}$$
($\Overline{a_i}$ denotes the residue of $a_i$).

\Item{(ii)}
$w=\{\chi,a_1,\dots,a_{r-2},\pi\}$ where $\chi$ is an element
of $H_p^1(K)$ corresponding to a cyclic extension
of $K$ of degree $p$ whose residue field is an inseparable extension of $k$
of degree $p$, $\pi$ is a prime element of $K$ 
and $a_1,\dots,a_{r-2}$ are elements of $\Cal O_K^*$ such that
$|k^p(\Overline{a_1},\dots,\Overline{a_{r-2}}):k^p|=p^{r-2}$.
\endRoster \enddf

\th Lemma 14

Let $K$ and $k$ be as in Definition 11.
Assume that $|k:k^p|=p^{r-1}$.
Then for every element $w\in H_p^r(K)\setminus H_{p,\ur}^r(K)$
there exists a finite extension $L$ of $K$ such that $|L:K|$ is prime to $p$
and such that the image of $w$ in $H_p^r(L)$ is standard.
\endth

\pf Proof

If $\chr(K)=p$ the proof goes just as in the proof of \cite{8, \S4 Lemma 5}
where the case of $r=2$ was treated.

If $\chr(K)=0$ we may assume that
$K$ contains a primitive $p$th root of 1.
Then the cohomological symbol
$h_{p,K}^r\colon K_r(K)/p\to H_p^r(K)$ is surjective
and
$$\text{\rm coker}(h_{p,K}^r\colon U_1K_r(K)\to H_p^r(K))\simeq \nu_{r-1}(k)$$
 by \cite{11} and $|k:k^p|=p^{r-1}$.

Here we are making the following:

\df Definition 12

Let $K$ be a complete discrete valuation field.
Then $U_iK_q(K)$ for $i,q\ge 1$ denotes
$U_i\Gamma_q(K)$ of Lemma 8 case (i) (take $A=\Cal O_K$ and
$B=K$).
\enddf


\df Definition 13

Let $k$ be a field of characteristic $p>0$.
As in Milne \cite{13}
denote by $\nu_{r}(k)$ the kernel of the homomorphism
$$\Omega_k^r\to \Omega_k^r/d(\Omega_k^{r-1}),
\quad
xdy_1/y_1\wedge\dots\wedge dy_r/y_r \mapsto (x^p-x)dy_1/y_1\wedge\dots\wedge dy_r/y_r.$$
\enddf

By \cite{11, Lemma 2} for every element $\alpha$ of $\nu_{r-1}(k)$
there is a finite extension $k'$ of $k$
such that 

\noindent $|k':k|$ is prime to $p$
and the image of $\alpha$ in   $\nu_{r-1}(k')$
is the sum of elements of type
$$dx_1/x_1\wedge\dots\wedge dx_r/x_r.$$
Hence we can follow the method of the proof of \cite{8, \S4 Lemma 5 
or \S2 Proposition 2}.
\qed\endpf

\pf Proof of Proposition 6

If $m$ is invertible in $k$ then
$H_m^r=H_{m,\ur}^r$.
Hence we may assume that
$\chr(k)=p>0$ and $m=p^i$, $i\ge 1$.
Since $\kr(p\colon H_{p^i}^r/H_{p^i,\ur}^r\to H_{p^i}^r/H_{p^i,\ur}^r)$
is isomorphic to
$H_{p}^r/H_{p,\ur}^r$ by \cite{11}, we may assume $m=p$.

The proof of part (1) is rather similar to the proof of part (2).
So we present here only the proof of part (2), but the method is directly
applicable to the proof of (1).

The proof is divided in several steps.
In the following $K$ always denotes an $n$-dimensional local field
with $n\ge 1$ and with residue field $k$
such that 

\noindent $\chr(k)=p>0$, except in Lemma 21.
$X$ denotes an object of $\Cal F_{\infty}$.

\smallskip

{\Em{Step 1}}. In this step $w$ denotes a standard element of $H_p^r(K)$
and $\Overline{w}$ is its image in $(H_p^r/H_{p,\ur}^r)(K)$.
We shall prove here that
$U_1\Gamma_q([X,\U K])\subset N(\Overline{w},[X,\U K])$.
We fix a presentation of $w$ as in (i) or (ii) of Definition 11.
Let $L$ be a cyclic extension of $K$ corresponding to $\chi$.
In case (i) (resp.\ (ii)) let $h$ be a prime element of $L$
(resp.\ an element of $\Cal O_L$ such that the residue class
$\Overline h$ is not contained in $k$).
Let $G$ be the subgroup of $K^*$ generated by $a_1,\dots,a_{r-1}$
(resp.\ by $a_1,\dots,a_{r-2},\pi$),  by $1+\Cal M_K$ and $N_{L/K}(h)$.
Let $l$ be the subfield of $k$ generated by the residue classes of $a_1,\dots,a_{r-1}$
(resp.\ $a_1,\dots,a_{r-2}, N_{L/K}(h)$).

Let $i\ge 1$. Let $G_{i,q}$ be the subgroup of $U_i\Gamma_q([X,\U K])$
generated by 

\noindent $\{U_i\Gamma_{q-1}([X,\U K]), G\}$ and $U_{i+1}\Gamma_q([X,\U K])$.
Under these notation we have the following Lemma 15, 16 ,17.

\th Lemma 15

\Roster 
\Item{(1)} $G_{i,q}\subset N_q(w,[X,\U K])+U_{i+1}\Gamma_q([X,\U K])$.

\Item{(2)} The homomorphism $\rho_i^q$ of Lemma 8 (3)
induces the surjections
$$[X,\U{\Omega_{k}^{q-1}}]\to [X,\U{\Omega_{k/l}^{q-1}}]@>\Overline{\rho_i^q}>>
U_1\Gamma_q([X,\U K])/G_{i,q}.$$

\Item{(3)} If $\rho_i^q$ is defined using a prime element $\pi$ which belongs to $G$
then the above homomorphism $\Overline{\rho_i^q}$
annihilates the image of the exterior derivation

\Item{} $d\colon [X,\U{\Omega_{k/l}^{q-2}}]\to [X,\U{\Omega_{k/l}^{q-1}}]$.
\endRoster 
\endth

\th Lemma 16

Let $a$ be an element of $K^*$ such that 
$v_K(a)=i$ and 

\noindent $a=a_1^{s(1)}\dots a_{r-1}^{s(r-1)}N_{L/K}(h)^{s(r)}$

\noindent (resp.\ $a=a_1^{s(1)}\dots a_{r-2}^{s(r-2)}\pi^{s(r-1)}N_{L/K}(h)^{s(r)}$)

\noindent where $s$ is a map $\{0,\dots,r\}\to\Bbb Z$ such that
$p\nmid s(j)$ for some $j\not=r$.

Then $1-x^pa\in N_1(w,[X,\U K])$ for each $x\in [X,\U{\Cal O_K}]$.
\endth
\pf Proof

It follows from the fact that
$w\in \{H_p^{r-1}(K),a\}$ and $1-x^pa$ is the norm
of $1-xa^{1/p}\in [X,\U{K(a^{1/p})}]^*$
($\U{K(a^{1/p})}$ denotes the ring object which represents the functor
$X\to K(a^{1/p})\otimes_K [X,\U K]$).
\qed \endpf

\th Lemma 17

Let $\sigma$ be a generator of $\Gal(L/K)$ and let
$a=h^{-1}\sigma(h)-1$, $b=N_{L/K}(a)$, $t=v_K(b)$.
Let $f=1$ in case (i) and let $f=p$ in case (ii).
Let $N\colon [X,\U L]^*\to  [X,\U K]^*$ be the norm homomorphism. Then:
\Roster 
\Item{(1)} If $f\vert i$ and $1\le i <t$ then for every
$x\in \Cal M_L^{i/f}[X,\U{\Cal O_L}]$ 
$$N(1+x)\equiv 1+ N(x) \mod \Cal M_K^{i+1}[X,\U{\Cal O_K}].$$

\Item{(2)} For every $x\in [X,\U{\Cal O_K}]$ 
$$N(1+xa)\equiv 1+ (x^p-x)b \mod \Cal M_K^{t+1}[X,\U{\Cal O_K}].$$
In case (ii) for every integer $r$ prime to $p$ and every
$x\in [X,\U{\Cal O_K}]$ 
$$N(1+xh^ra)\equiv 1+ x^pN(h)^rb \mod \Cal M_K^{t+1}[X,\U{\Cal O_K}].$$

\Item{(3)} 
$$1+\Cal M_K^{t+1}[X,\U{\Cal O_K}]\subset N(1+\Cal M_L^{t/f+1}[X,\U{\Cal O_L}]).$$
\endRoster 
\endth

\pf Proof

Follows from the computation of the norm homomorphism
$L^*\to K^*$ in Serre \cite{15, Ch. V \S3}
and \cite{8, \S1}.
\qed\endpf

From these lemmas we have
\Roster 
\Item{(1)} If $0<i<t$ then
$$U_i\Gamma_q([X,\U K])\subset N_q(w,[X,\U K])+U_{i+1}\Gamma_q([X,\U K]).$$

\Item{(2)} $U_{t+1}\Gamma_q([X,\U K])\subset N_q(w,[X,\U K])$. 

\Item{(3)} In case (ii) let $a_{r-1}=N_{L/K}(h)$.
then in both cases (i) and (ii)
the homomorphism
$$
\aligned
&[X,\U{\Omega_k^{q+r-2}}]\to U_1\Gamma_q([X,\U K])/N_q(w,[X,\U K]),\\
&xd\Overline{a_1}/\Overline{a_1}\wedge\dots\wedge 
d\Overline{a_{r-1}}/\Overline{a_{r-1}}\wedge dy_1/y_1\wedge\dots\wedge
dy_{q-1}/y_{q-1}\\
&\qquad\qquad\qquad\qquad \mapsto \{1+\wt{x}b,\wt{y_1},\dots,\wt{y_{q-1}}\},
\endaligned
$$
($x\in [X,\U k], y_i\in [X,\U{k^*}]$) 
annihilates $(1-\Car )[X,\U{\Omega_k^{q+r-2}}]$.
\endRoster

 Lemma 7 and (1), (2), (3) imply that
$U_1\Gamma_q([X,\U K])$ is contained
in the sum of $N_q(w,[X,\U K])$ and the image of
$\lcf(X,U_{t+1}K_q(K))$.

\th Lemma 18

For each $u\in \Cal O_K$ there exists an element
$\psi$ of $H_{p,\ur}^1(K)$ such that
$(1+ub)N_{L/K}(h)^{-1}$ is contained
in the norm group $N_{L'/K}{L'}^*$ where
$L'$ is the cyclic extension of $K$ corresponding to $\chi+\psi$
{\rm(}$\chi$ corresponds to $L/K${\rm)}.
\endth

\pf Proof

Follows from \cite{9, \S3.3 Lemma 15}
(can be proved using the formula
$$N_{L_{\ur}/K_{\ur}}(1+xa)\equiv 1+(x^p-x)b \mod b\Cal M_{K_{\ur}}$$
for $x\in \Cal O_{K_{\ur}}$.
\qed\endpf

Lemma 18 shows that $1+ub$ is contained in the subgroup generated by
$N_{L/K}L^*$ and $N_{L'/K}{L'}^*$, $\chi_L=0$, $\chi_{L'}\in H_{p,\ur}^1(L')$.

\smallskip

{\Em{Step 2}}. Next we prove that
$$U_1\Gamma_q([X,\U K])\subset N(\Overline{w},[X,\U K])$$
 for every $w\in H_p^r(K)$ where $\Overline{w}$ is the image of $w$ in
$(H_p^r/H_{p,\ur}^r)(K)$.

\th Lemma 19

Let $q,r\ge 1$ and let $w\in H_p^r(K)$. Then
there exists $i\ge 1$ such that
$p^i\Gamma_q([X,\U{K'}])$ and
$U_{e(K'|K)i}\Gamma_q([X,\U{K'}])$
are contained in 
 $N_q(w_{K'},[X,\U{K'}])$ for every $K'\in \Cal E(K)$
where $w_{K'}$ denotes the image of $w$ in $H_p^r(K')$ and $e(K'|K)$ denotes
the ramification index of $K'/K$.
\endth

\th Lemma 20

Let $i\ge 1$ and $x\in U_1\Gamma_q([X,\U{K}])$; 
(resp.\ $x=\{u_1,\dots,u_q\}$ with 

\noindent $u_i\in [X,\U{\Cal O_K^*}]$;
resp.\ $x\in \Gamma_q([X,\U{K}])$).

Then there exists a triple $(J,0,f)$ and a family
$(E_j,x_j)_{j\in J}$ which satisfy all the conditions of Definition 6
except
(iv) and satisfy condition (iv)' below.
\Roster 
\Item{(iv)'} If $j\not\in f(J)$ then $x_j$ satisfy one of the following three properties: 
\ItemItem{(a)} $x_j\in p^i\Gamma_q([X,\U{E_j}])$.

\ItemItem{(b)} $x_j\in U_{e(E_j|K)i}\Gamma_q([X,\U{E_j}])$;
(resp.\ (b) $x_j\in U_1\Gamma_q([X,\U{E_j}])$.

\ItemItem{(c)} Let $\Overline{E_j}$ be the residue field of $E_j$.
There are elements $c_1,\dots,c_{q-1}$ of $\Cal O_{E_j}^*$
such that

\ItemItem{} $x_j\in \{U_1\Gamma_1([X,\U{E_j}]),c_1,\dots,c_{q-1}\}$
and 
$|\Overline{E_j}^p(c_1,\dots,c_{q-1}):\Overline{E_j}^p|=p^{q-1}$;

\ItemItem{}(resp.\ (c) There are elements $b_1,\dots,b_q$ of $[X,\U{\Cal O_{E_j}^*}]$ and 
$c_1,\dots,c_{q}$ of $\Cal O_{E_j}^*$
such that  
 $x_j=\{b_1,\dots,b_q\}$ and such that for each $m$ the residue class
$\Overline{b_m}\in [X,\U{\Overline{E_j}}]$ belongs to
$[X,\U{\Overline{E_j}}]^p[\Overline{c_m}]$ and 
$|\Overline{E_j}^p(c_1,\dots,c_{q}):\Overline{E_j}^p|=p^{q}$);

\ItemItem{}{\rm(}resp.\ (c) There are elements $b_1,\dots,b_{q-1}$ of $[X,\U{\Cal O_{E_j}^*}]$ and 
$c_1,\dots,c_{q-1}$ of $\Cal O_{E_j}^*$
such that
$x_j\in \{[X,\U{E_j}]^*, b_1,\dots,b_{q-1}\}$
and such that for each $m$ the residue class
$\Overline{b_m}\in [X,\U{\Overline{E_j}}]$ belongs to
$[X,\U{\Overline{E_j}}]^p[\Overline{c_m}]$ and \hfill\break
$|\Overline{E_j}^p(c_1,\dots,c_{q-1}):\Overline{E_j}^p|=p^{q-1}${\rm)}.
\endRoster\endth

Using Lemma 19 and 20
it suffices for the purpose of this step to consider
the following elements

\noindent $\{u,c_1,\dots,c_{q-1}\}\in U_1\Gamma_q([X,\U K])$ such that 
$u\in U_1\Gamma_1([X,\U K])$, $c_1,\dots,c_{q-1}\in \Cal O_K^*$
and $|k^p(\Overline{c_1},\dots, \Overline{c_{q-1}}:k^p|=p^{q-1}$.

For each $i=1,\dots,q-1$ and each $s\ge 0$ take a $p^s$th
root $c_{i,s}$ of $-c_i$
satisfying $c_{i,s+1}^p=c_{i,s}$.
Note that
$N_{k(c_{i,s+1})/k(c_{i,s})}(-c_{i,s+1})=-c_{i,s}$.
For each $m\ge 0$ write $m$ in the form
$(q-1)s+r$ ($s\ge 0$, $0\le r<q-1$).
Let $E_m$ be the finite extension of $K$ of degree $p^m$ generated by
$c_{i,s+1}$ ($1\le i\le r$) and $c_{i,s}$ ($r+1\le i\le q-1$) and let
$$x_m=\{u,-c_{1,s+1},\dots, -c_{r,s+1},-c_{r+1,s},\dots,-c_{q-1,s}\}\in \Gamma_q([X,\U{E_m}]).
$$
Then  $E_\infty=\inlim E_m$ is a henselian discrete valuation field
with residue field $\Overline{E_\infty}$
satisfying 
$|\Overline{E_\infty}:\Overline{E_\infty}^p|\le p^{r-1}$. 
Hence by Lemma 14 and Lemma 21 below
there exists $m<\infty$ such that
for some finite extension
$E_m'$ of $E_m$ of degree prime to $p$
the image of $w$ in $H_p^r(E_m')$ is standard.
Let $J=\{0,1,\dots,m,m'\}$, $f(j)=j-1$ for $1\le j\le m$,
$f(0)=0$, $f(m')=m$, $E_{m'}=E_m'$ and
$$x_{m'}=\{u^{1/|E_m':E_m|}, c_1,\dots,c_{q-1}\}.$$
Then from Step 1 we deduce
$\{u,c_1,\dots,c_{q-1}\}\in N_q(\Overline{w},[X,\U K])$.

\th Lemma 21

Let $K$ be a henselian discrete valuation field,
and let $\wh{K}$ be its completion.
Then $H_m^r(K)\simeq H_m^r(\wh{K})$ for every $r$ and $m$.
\endth

\pf Proof

If $m$ is invertible in $K$ this follows from the isomorphism
$\Gal(\wh{K}^{\sep}/\wh{K})\simeq \Gal(K^{\sep}/K)$
(cf.\ \cite{1, Lemma 2.2.1}).
Assume $\chr(K)=p>0$ and $m=p^i$ ($i\ge 1$).
For a field $k$ of characteristic $p>0$
the group $H_{p^i}^r(k)$ is isomorphic to
$(H_{p^i}^1(k)\otimes k^*\otimes\dots\otimes k^*)/J$ where
$J$ is the subgroup of the tensor product
generated by elements
of the form (cf.\ \cite{9, \S2.2 Corollary 4 to Proposition 2})
\Roster 
\Item{(i)} $\chi\otimes a_1\otimes\dots\otimes a_{r-1}$ such that $a_i=a_j$ for some $i\not=j$,

\Item{(ii)} $\chi\otimes a_1\otimes\dots\otimes a_{r-1}$ such that 
$a_i\in N_{k_\chi/k}k_\chi^*$
for some $i$ where $k_\chi$ is the extension of $k$ corresponding to $\chi$.
\endRoster

By the above isomorphism of the Galois groups
$H_{p^i}^1(K)\simeq  H_{p^i}^1(\wh{K})$.
Furthermore if $L$ is a cyclic extension of $K$
then
$1+\Cal M_K^n\subset N_{L/K}L^*$
and $1+\Cal M_{\wh{K}}^n\subset N_{L\wh{K}/\wh{K}}(L\wh{K})^*$
for sufficiently large $n$.
Since
$K^*/(1+\Cal M_K^n)\simeq \wh{K}^*/(1+\Cal M_{\wh{K}}^n)$,
the lemma follows.
\qed\endpf

{\Em{Step 3}}.
In this step we prove that the subgroup of $\Gamma_q([X,\U K])$ generated by

\noindent $U_1\Gamma_q([X,\U K])$ and elements of the form
$\{u_1,\dots,u_q\}$
($u_i\in [X,\U{\Cal O_K^*}]$)
is contained in $N_q(\Overline{w},[X,\U K])$.
By Lemma 20 it suffices to consider
elements $\{b_1,\dots,b_q\}$ such that
$b_i\in [X,\U{\Cal O_K^*}]$ and such that
there are elements
$c_j\in \Cal O_K^*$
satisfying
$$|k^p(\Overline{c_1},\dots,\Overline{c_q}):k^p|=p^q$$ and 
$\Overline{b_i}\in [X,\U k]^p[\Overline{c_i}]$ for each $i$.
Define fields $E_m$ as in Step 2 replacing $q-1$ by $q$.
Then
$E_\infty=\inlim E_m$ is a henselian discrete valuation field
with residue field $\Overline{E_\infty}$ satisfying
$|\Overline{E_\infty}:\Overline{E_\infty}^p|\le p^{r-2}$.
Hence
$H_p^r(\wh{E_\infty})=H_{p,\ur}^r(\wh{E_\infty})$.
By Lemma 21 there exists $m<\infty$
such that $w_{E_m}\in H_{p,\ur}^r(E_m)$.

\smallskip

{\Em{Step 4}}.
Let $w$ be a standard element.
Then there exists a prime element $\pi$
of $K$ such that
$\pi\in N_1(w,[X,\U K])=\Gamma_q([X,\U K])$.

\smallskip

{\Em{Step 5}}.
Let $w$ be any element of $H_p^r(K)$.
To show that
$\Gamma_q([X,\U K])=N_q(\Overline{w},[X,\U K])$ it suffices using Lemma 20
to consider elements of $\Gamma_q([X,\U K])$
of the form

\noindent $\{x,b_1,\dots,b_{q-1}\}$ ($x\in [X,\U K]^*$, $b_i\in [X,\U{\Cal O_K^*}]$)
such that 
there are elements

\noindent $c_1,\dots,c_{q-1}\in \Cal O_K^*$ satisfying  
$|k^p(\Overline{c_1},\dots,\Overline{c_{q-1}}):k^p|=p^{q-1}$ and 
$\Overline{b_i}\in [X,\U k]^p[\Overline{c_i}]$ for each $i$.
The fields $E_m$ are defined again as in Step 2,
and we are reduced to the case where $w$ is standard.
\qed\endth

\HH 5. Proof of Theorem 2

\phantom{}
\par

Let $K$ be an $n$-dimensional local field.
By \cite{9, \S3 Proposition 1}
$H^r(K)=0$ for 

\noindent $r>n+1$
and there exists a canonical isomorphism
$H^{n+1}(K)\simeq \Bbb Q/\Bbb Z$.

\noindent For $0\le r\le n+1$ the canonical pairing
$$\{\,\,,\,\,\}\colon H^r(K)\times K_{n+1-r}(K)\to H^{n+1}(K)$$
(see subsection 4.2)
induces a homomorphism
$$\Phi_K^r\colon H^r(K)\to \Hom(K_{n+1-r}(K),\Bbb Q/\Bbb Z).$$
if $w\in H^r(K)$ with $r\ge 1$ (resp.\ $r=0$)
then $\Phi_K^r(w)$ annihilates the norm group
$N_{n+1-r}(w)$ (resp.\ $\Phi_K^r(w)$ annihilates $mK_{n+1}(K)$ where $m$
is the order of $w$).
Since $N_{n+1-r}(w)$ (resp.\ $mK_{n+1}(K)$)
is open in $K_{n+1-r}(K)$ by Theorem 3 (resp.\ Corollary to Lemma 12),
$\Phi_K^r(w)$ is a continuous character of $K_{n+1-r}(K)$ of finite order.

\HHH 5.1. Continuous characters of prime order

\phantom{}\smallskip\par

In this subsection we shall prove that for every prime $p$
the map $\Phi_K^r$ ($0\le r\le n+1$)
induces a bijection between $H_p^r(K)$ (cf.\ Definition 10)
and the group of all continuous characters of order $p$ of $K_{n+1-r}(K)$.
We may assume that $n\ge 1$ and $1\le r\le n$.
Let $k$ be the residue field of $K$.
In the case where $\chr(k)\not=p$ the above assertion
follows by induction on $n$ from the isomorphisms
$$
H_p^r(k)\oplus H_p^{r-1}(k)\simeq H_p^r(K),\quad
K_q(k)/p\oplus K_q(k)/p\simeq K_q(K)/p.
$$
Now we consider the case of $\chr(k)=p$.

\df Definition 14

Let $K$ be a complete discrete valuation field with residue field $k$ of characteristic $p>0$.
For $r\ge 1$ and $i\ge 0$ we define the subgroup
$T_iH_p^r(K)$ of $H_p^r(K)$ as follows.
\Roster 
\Item{(1)} If $\chr(K)=p$ then let
$\delta_K^r\colon \Omega_K^{r-1}=C_1^{r-1}(K)\to H_p^r(K)$ be the canonical
projection.
Then $T_iH_p^r(K)$ is the subgroup of $H_p^r(K)$ generated by elements
of the form
$$
\delta_K^r(xdy_1/y_1\wedge\dots\wedge dy_{r-1}/y_{r-1}),
\quad x\in K, y_1,\dots, y_{r-1}\in K^*, v_K(x)\ge -i.
$$

\Item{(2)} If $\chr(K)=0$ then let
$\zeta$ be a primitive $p$th root of 1,
and let $L=K(\zeta)$.

\Item{} Let $j=(pe_K/(p-1)-i)e(L|K)$ where
$e_K=v_K(p)$ and $e(L|K)$ is the ramification index of $L/K$.
If $j\ge 1$ let
$U_jH_p^r(L)$ be the image of $U_j K_r(L)$
(cf.\ Definition 12)
under the cohomological symbol
$K_r(L)/p\to H_p^r(L)$.
If $j\le 0$, let $U_jH_p^r(L)=H_p^r(L)$.
Then $T_iH_p^r(K)$ is the inverse image of $U_j H_p^r(L)$ under the canonical injection
$H_p^r(K)\to H_p^r(L)$.
\endRoster 
\enddf

\rk Remark

$T_iH_p^1(K)$ coincides with the subgroup consisting
of elements which corresponds to cyclic extensions of $K$ of degree $p$
with ramification number $\le i$
(the ramification number is defined as $t$ of Lemma 17).
\endrk

Let $K$ be as in Definition 14, and assume that $|k:k^p|<\infty$.
Fix $q,r\ge 1$ such that
$|k:k^p|=p^{q+r-2}$.
Let $T_i=T_iH_p^r(K)$, for $i\ge0$;
let $U_i$ be the image of $U_iK_q(K)$ in $K_q(K)/p$ for $i\ge 1$,
and let $U_0=K_q(K)/p$.
Let $e=v_K(p)$ ($=\infty$ if
$\chr(K)=p$).
Fix a prime element $\pi$ of $K$.
Via the homomorphism
$$(x,y)\mapsto \rho_i^q(x)+\{\rho_i^{q-1}(y),\pi\}$$
of Lemma 8 whose kernel is known by \cite{11},
we identify $U_i/U_{i+1}$ with the following groups:
\Roster 
\Item{(1)} $K_q(k)/p\oplus K_{q-1}(k)/p$ if $i=0$.

\Item{(2)} $\Omega_k^{q-1}$ if $0<i<pe/(p-1)$ and $i$ is prime to $p$.

\Item{(3)} $\Omega_k^{q-1}/\Omega_{k,d=0}^{q-1}\oplus \Omega_k^{q-2}/\Omega_{k,d=0}^{q-2}$
if $0<i<pe/(p-1)$ and $p|i$.

\Item{(4)} $\Omega_k^{q-1}/D_{a,k}^{q-1}\oplus \Omega_k^{q-2}/D_{a,k}^{q-2}$ if
$\chr(K)=0$, $pe/(p-1)$ is an integer and $i=pe/(p-1)$.

\Item{(5)} 0 if $i>pe/(p-1)$.
\endRoster 

Here in (3) $\Omega_{k,d=0}^q$ ($q\ge0$)
denotes the kernel of the exterior derivation

\noindent $d\colon \Omega_k^{q}\to \Omega_k^{q+1}$.
In (4) $a$ denotes the residue class of $p\pi^{-e}$
where $e=v_K(p)$ and
$D_{a,k}$ denotes the subgroup of $\Omega_k^q$ generated by
$d(\Omega_k^{q-1})$ and elements
of the form
$$(x^p+ax)dy_1/y_1\wedge\dots\wedge dy_q/y_q.$$

Note that $H_p^{r+1}(K)\simeq H_p^{q+r-1}(k)$ by \cite{11}.
Let $\delta=\delta_k^{q+r-1}\colon \Omega_k^{q+r-2}\to H_p^{q+r-1}(k)$
(Definition 14).

\th Lemma 22

In the canonical pairing
$$H_p^r(K)\times K_q(K)/p\to H_p^{q+r}(K)\simeq H_p^{q+r-1}(k)$$
$T_i$ annihilates $U_{i+1}$ for each $i\ge 0$.
Furthermore,
\Roster 
\Item{(1)} $T_0=H_{p,\ur}^r(k)\simeq H_p^r(k)\oplus H_p^{r-1}(k)$, 
and the induced pairing
$$T_0\times U_0/U_1\to H_p^{q+r-1}(k)$$
is identified with the direct sum of the canonical
pairings
$$
H_p^r(k)\times K_{q-1}(k)/p\to H_p^{q+r-1}(k),\quad 
H_p^{r-1}(k)\times K_{q}(k)/p\to H_p^{q+r-1}(k).$$

\Item{(2)} If $0<i<pe/(p-1)$ and $i$ is prime to $p$ then there exists
an isomorphism
$$T_i/T_{i-1}\simeq \Omega_k^{r-1}$$ such that
the induced pairing
$T_i/T_{i-1}\times U_i/U_{i+1}\to H_p^{q+r-1}(k)$
is identified with
$$\Omega_k^{r-1}\times \Omega_k^{q-1}\to H_p^{q+r-1}(k), \quad (w,v)\mapsto \delta(w\wedge v).$$

\Item{(3)} If $0<i<pe/(p-1)$ and $p|i$ then there exists an isomorphism
$$T_i/T_{i-1}\simeq \Omega_k^{r-1}/\Omega_{k,d=0}^{r-1}\oplus
\Omega_k^{r-2}/\Omega_{k,d=0}^{r-2}$$
such that the induced pairing is identified with 
$$(w_1\oplus w_2,v_1\oplus v_2)\mapsto \delta(dw_1\wedge v_2+dw_2\wedge v_1).$$

\Item{(4)} If $\chr(K)=0$ and $pe/(p-1)$ is not an integer,
then $H_p^r(K)=T_i$ for the maximal integer $i$ smaller than $pe/(p-1)$.
Assume that $\chr(K)=0$ and $pe/(p-1)$ is an integer.
Let $a$ be the residue element of $p\pi^{-e}$ and let for $s\ge 0$
$$\nu_s(a,F)=\kr(\Omega_{k,d=0}^s\to \Omega_k^s,\quad w\mapsto 
\text{\tenrm C}(w)+aw)$$
{\rm(}$\text{\tenrm C}$ denotes the Cartier operator{\rm)}.
Then there exists an isomorphism
$$T_{pe/(p-1)}/T_{pe/(p-1)-1}\simeq \nu_r(a,k)\oplus \nu_{r-1}(a,k)$$
such that the induced pairing
is identified with
$$(w_1\oplus w_2,v_1\oplus v_2)\mapsto \delta(w_1\wedge v_2+w_2\wedge v_1).$$
\endRoster
\endth
\pf Proof

If $\chr(K)=p$ the lemma follows from a computation in the differential modules $\Omega_K^s$ ($s=r-1, q+r-1$).
In the case where $\chr(K)=0$ let
$\zeta$ be a primitive $p$th root of 1
and let $L=K(\zeta)$.
Then the cohomological symbol
$K_r(L)/p\to H_p^r(L)$ is surjective
and the structure
of $H_p^r(L)$ is explicitly
given in \cite{11}.
Since
$$H_p^r(K)\simeq \{x\in H_p^r(L): \sigma(x)=x\quad\text{\rm for all $\sigma\in \Gal(L/K)$}\},$$
the structure of $H_p^r(K)$ is deduced from that of $H_p^r(L)$ and the description of the pairing
$$H_p^r(K)\times K_q(K)/p\to H_p^{q+r}(K)$$
follows from a computation of the pairing 
$$K_r(L)/p\times K_q(L)/p\to K_{q+r}(L)/p.\eqno{\qed}$$
\endpf

\th Lemma 23

Let $K$ be an $n$-dimensional local field
such that $\chr(K)=p>0$.
Then the canonical map
$\delta_K^n\colon \Omega_K^n\to H_p^{n+1}(K)\simeq \Bbb Z/p$
{\rm(}cf.\ Definition 14{\rm)}
comes from a morphism $\U{\Omega_K^n}\to \Bbb Z/p$ of $\Cal A_\infty$.
\endth
\pf Proof

Indeed it comes from the composite morphism of $\Cal F_\infty$
$$\U{\Omega_K^n}@>\theta_2>> k_0@> \Tr_{k_0/\Bbb F_p}>> \Bbb F_p$$
defined by Lemma 7.
\qed\endpf

Now let $K$ be an $n$-dimensional local field ($n\ge 1$)
with residue  field $k$ such that
$\chr(k)=p>0$.
Let $1\le r\le n$, $q=n+1-r$,
and let $T_i$ and $U_i$ ($i\ge 0$)
be as in Lemma 22.

The injectivity of the map induced by
$\Phi_K^r$
$$H_p^r(K)\to \Hom(K_{n+1-r}(K)/p, \Bbb Z/p)$$
follows by induction on $n$
from the injectivity of
$T_i/T_{i-1}\to \Hom(U_i/U_{i+1},\Bbb Z/p)$, $i\ge 1$.
Note that this injectivity for all prime $p$
implies the injectivity of $\Phi_K^r$.

Now let $\varphi\colon K_{n+1-r}(K)\to \Bbb Z/p$
be a continuous character of order $p$.
We prove that there is an element
$w$ of $H_p^r(K)$ such that
$\varphi=\Phi_K^r(w)$.

The continuity of $\varphi$ implies that
there exists $i\ge 1$ such that
$$\varphi(\{x_1,\dots,x_{n+1-r}\})=0\quad \text{\rm for all $x_1,\dots,x_{n+1-r}\in 1+\Cal M_K^i$.}$$
Using Graham's method \cite{6, Lemma 6}
we deduce that
$\varphi(U_i)=0$ for
some $i\ge 1$.
We prove the following assertion ($A_i$) ($i\ge 0$)
by downward induction on $i$.

\noindent ($A_i$) The restriction of $\varphi$ to $U_i$
coincides with the restriction of
$\Phi_K^r(w)$ for some $w\in H_p^r(K)$.

Indeed, by induction on $i$ there exists $w\in H_p^r(K)$
such that
the continuous character
$\varphi'=\varphi-\Phi_K^r(w)$ annihilates $U_{i+1}$.

In the case where $i\ge 1$ the continuity of $\varphi'$ implies that the map
$$\Omega_k^{n-r}\oplus \Omega_k^{n-r-1}@>\text{\rm Lemma 8}>>
U_i/U_{i+1}@>\varphi'>>\Bbb Z/p$$
comes from a morphism of $\Cal F_\infty$.
By additive duality of Proposition 3 and Lemma 23 applied to $k$
the above map is expressed in the form
$$(v_1,v_2)\mapsto \delta_k^n(w_1\wedge v_2+w_2\wedge v_1)$$
for some $w_1\in \Omega_k^n, w_2\in \Omega_k^{r-1}$.
By the following argument
the restriction of $\varphi'$ to $U_i/U_{i+1}$
is induced by an element
of $T_i/T_{i-1}$.
For example, assume $\chr(K)=0$ and $i=pe/(p-1)$
(the other cases are treated similarly and more easily).
Since $\varphi'$ annihilates
$d(\Omega_k^{n-r-1})\oplus d(\Omega_k^{n-r-2})$
and $\delta_k^n$ annihilates $d(\Omega_k^{n-2})$ we get
$$\delta_k^n(dw_1\wedge v_2)=\pm \delta_k^n(w_2\wedge dv_2)=0\quad
\text{\rm for all $v_2$.}$$
Therefore $dw_1=0$.
For every $x\in F$, $y_1,\dots,y_{n-r-1}\in F^*$ we have
$$\delta_k^n\bigl((\Car (w_1)+aw_1)\wedge x\frac{dy_1}{y_1}\wedge\dots\wedge \frac{dy_{n-r-1}}{y_{n-r-1}}\bigr)=
\delta_k^n\bigl(w_1\wedge (x^p+ax)\frac{dy_1}{y_1}\wedge\dots\wedge \frac{dy_{n-r-1}}{y_{n-r-1}}\bigr)=0$$
(where $a$ is as in Lemma 22 (4)).
Hence $w_1\in \nu_r(a,k)$ and similarly $w_2\in \nu_{r-1}(a,k)$.

In the case where $i=0$ Lemma 22 (1) and induction on $n$ imply that
there is an element $w\in T_0$ such that
$\varphi'=\Phi_K^r(w)$.

\HHH 5.2. Continuous characters of higher orders

\phantom{}\smallskip\par

In treatment of continuous characters of higher order
the following proposition will play a key role.

\th Proposition 7

Let $K$ be an $n$-dimensional local field.
Let $p$ be a prime number distinct from the characteristic of $K$.
Assume that $K$ contains a primitive $p$th root $\zeta$ of 1.
Let $r\ge0$ and $w\in H^r(K)$.
Then the following two conditions are equivalent.
\Roster 
\Item{(1)} $w=pw'$ for some $w'\in H^r(K)$.

\Item{(2)} $\{w,\zeta\}=0$ in $H^{r+1}(K)$.
\endRoster\endth
\pf Proof

We may assume that $0\le r\le n$.
Let $\delta_r\colon H^r(K)\to H^{r+1}(K,\Bbb Z/p)$
be the connecting homomorphism induced by the exact sequence of
$\Gal(K^{\sep}/K)$-modules
$$0\to \Bbb Z/p\to \inlim_i \, \mu_{p^i}^{\otimes(r-1)}@>p>>
\inlim_i\, \mu_{p^i}^{\otimes(r-1)}\to 0.$$
Condition (1) is clearly equivalent to $\delta_r(w)=0$.

First we prove the proposition
in  the case where $r=n$.
Since the kernel of 
$$\delta_n\colon H^n(K)\to H^{n+1}(K,\Bbb Z/p)\simeq \Bbb Z/p$$
is contained in the kernel of the homomorphism
$\{\,\,\,,\zeta\}\colon H^n(K)\to H^{n+1}(K)$
it suffices to prove that the latter homomorphism is not a zero map.
Let $i$ be the maximal natural number
such that $K$ contains a primitive $p^i$th root of 1.
Since the image $\chi$ of a primitive $p^i$th root of 1 under
the composite map 
$$K^*/K^{*p}\simeq H^1(K,\mu_p)\simeq H^1(K,\Bbb Z/p)\to H^1(K)$$
is not zero, the injectivity of $\Phi_K^1$ shows that
there is an element
$a$ of $K_n(K)$ such that
$\{\chi,a\}\not=0$.
Let $w$ be the image of $a$ under the composite map
induced by the cohomological symbol
$$K_n(K)/p^i\to H^n(K,\mu_{p^i}^{\otimes n})\simeq H^n(K,\mu_{p^i}^{\otimes(n-1)})\to H^n(K).$$
Then $\{\chi,a\}=\pm\{w,\zeta\}$.

Next we consider the general case of $0\le r\le n$.
Let $w$ be an element of $H^r(K)$ such that
$\{w,\zeta\}=0$.
Since the proposition holds for $r=n$
we get $\{\delta_r(w),a\}=\delta_n(\{w,a\})=0$ for all $a\in K_{n-r}(K)$.
The injectivity of $\Phi_K^{r+1}$ implies
$\delta_r(w)=0$.
\qed\endpf

\rk Remark

We conjecture that condition (1) is equivalent to condition (2)
for every field $K$.

This conjecture is true if $\oplus_{r\ge 1}H^r(K)$ is generated by 
$H^1(K)$ as  

\noindent a $\oplus_{q\ge 0}K_q(K)$-module.
\endrk

\pf Completion of the proof of Theorem 2

Let $\varphi$ be a non-zero continuous character of $K_{n+1-r}(K)$ of finite order,
and let $p$ be a prime divisor of the order of $\varphi$.
By induction on the order there exists
an element $w$ of $H^r(K)$
such that
$p\varphi=\Phi_K^r(w)$.
If $\chr(K)=p$ then $H^r(K)$ is $p$-divisible.
If $\chr(K)\not=p$, let $L=K(\zeta)$ where $\zeta$
is a primitive $p$th root of 1 and let
$w_L$ be the image of $w$ in $H^r(L)$.
Then $\Phi_L^r(w_L)\colon K_{n+1-r}(L)\to \Bbb Q/\Bbb Z$
coincides with the composite
$$K_{n+1-r}(L)@>N_{L/K}>> K_{n+1-r}(K)@>p\varphi>> \Bbb Q/\Bbb Z$$
and hence
$\{w_L,\zeta,a\}=0$ in $H^{n+1}(L)$ for all $a\in K_{n-r}(L)$.
The injectivity of $\Phi_L^{r+1}$ and Proposition 7 imply that
$w_L\in pH^r(L)$.
Since $|L:K|$ is prime to $p$, $w$ belongs to $pH^r(K)$.

Thus there is an element $w'$ of $H^r(K)$ such that $w=pw'$.
Then $\varphi-\Phi_K^r(w')$ is a continuous character annihilated by $p$.
\qed\endpf

\medskip

 \Bib        References

 \rf {1} M. Artin, Dimension cohomologique; premiers r\'esultats,
Th\'eorie des topos et cohomologie etale des sch\'emas, Tome 3, Expos\'e X, Lecture Notes in Math. 305, Springer, Berlin 1973, 43--63

\rf{2} H. Bass, J. Tate,
The Milnor ring of a global field,
Algebraic $K$-theory II, Lecture Notes in Math. 342, Springer, Berlin 1972,
349--446.

\rf{3} S. Bloch,
Algebraic $K$-theory and crystalline cohomology,
Publ. Math. IHES 47, 1977, 187--268.

\rf{4} P. Cartier, Une nouvelle op\'eration sur les formes diff\'erentielles,
C. R. Acad. Sc. Paris 244, 1957, 426--428.

\rf{5} P. Deligne, Cohomologie \`a support propre et construction du functeur
$f^{!}$,
\ in \ R. Har\-ts\-horne, Residue and duality, Lecture Notes in Math. 20,
Springer, Berlin 1966, 404--421.

\rf{6} J. Graham,
Continuous symbols on fields of formal power series,
Algebraic $K$-theory II, Lecture Notes in Math. 342, Springer, Berlin 1972,
474--486.

\rf{7} A. Grothendieck,
Elements de G\'eom\'etrie Alg\'ebrique IV,
Premi\`ere Partie, Publ. Math. IHES 20, 1964.

\rf{8} K. Kato, A generalization of local class field theory
by using $K$-groups I,
J. Fac. Sci. Univ. Tokyo Sec. IA 26 No.2, 1979, 303--376.

\rf{9} K. Kato, A generalization of local class field theory
by using $K$-groups II,
J. Fac. Sci. Univ. Tokyo Sec. IA 27 No.3, 1980, 603--683.

\rf{10} K. Kato, A generalization of local class field theory
by using $K$-groups III,
J. Fac. Sci. Univ. Tokyo Sec. IA 29 No.1, 1982, 31--42.

\rf{11} K. Kato, Galois cohomology of complete discrete valuation fields,
Algebraic $K$-theory, Part II (Oberwolfach, 1980), Lecture Notes in Math. 967,  1982, 215--238.

\rf{12} S. Lefschetz, Algebraic Topology, Amer. Math. Soc. Colloq. Publ., 1942.

\rf{13} J. S. Milne, Duality in flat cohomology of a surface,
Ann. Sc. Ec. Norm. Sup.
4\`eme s\'erie 9, 1976, 171--202.

\rf{14} J. Milnor,
Algebraic $K$-theory and quadratic forms,
Invent. Math. 9, 1970, 318--344.

\rf{15} J.-P. Serre, Corps Locaux, Hermann, Paris 1962.

\rf{16} J.-P. Serre, Cohomologie Galoisienne,
Lecture Notes in Math. 5, Springer, Berlin 1965.

\rf{17} A. Weil, Basic Number Theory, Springer, Berlin 1967.

 \endBib
 
 \Coordinates

Department of Mathematics \ 
University of Tokyo

3-8-1 Komaba Meguro-Ku Tokyo 153-8914 Japan
\endCoordinates

\vfill
\pagebreak

\end

%% file: m3-macs.tex
\expandafter\ifx\csname mthreemacsloaded\endcsname\relax\else \fi

\magnification1100
\input amstex


 \catcode`\@=11
 \let\wlog@ld\wlog
 \def\wlog#1{\relax}

 \newif\ifIN@
 \def\m@rker{\m@@rker}
 \def\IN@{\expandafter\INN@\expandafter}
 \long\def\INN@0#1@#2@{\long\def\NI@##1#1##2##3\ENDNI@
    {\ifx\m@rker##2\IN@false\else\IN@true\fi}%
     \expandafter\NI@#2@@#1\m@rker\ENDNI@}
  \newtoks\Initialtoks@  \newtoks\Terminaltoks@
  \def\SPLIT@{\expandafter\SPLITT@\expandafter}
  \def\SPLITT@0#1@#2@{\def\TTILPS@##1#1##2@{%
     \Initialtoks@{##1}\Terminaltoks@{##2}}\expandafter\TTILPS@#2@}
  \newtoks\Trimtoks@

 \def\ForeTrim@{\expandafter\ForeTrim@@\expandafter}
 \def\ForePrim@0 #1@{\Trimtoks@{#1}}
 \def\ForeTrim@@0#1@{\IN@0\m@rker. @\m@rker.#1@%
     \ifIN@\ForePrim@0#1@%
     \else\Trimtoks@\expandafter{#1}\fi}
 
  \def\Trim@0#1@{%
      \ForeTrim@0#1@%
      \IN@0 @\the\Trimtoks@ @%
        \ifIN@
             \SPLIT@0 @\the\Trimtoks@ @\Trimtoks@\Initialtoks@
             \IN@0\the\Terminaltoks@ @ @%
                 \ifIN@
                 \else \Trimtoks@ {FigNameWithSpace}%
                 \fi
        \fi
      }

  \font\titlebold=cmbx12 scaled 1200
  \font\twelvebold=cmbx12
  \font\tenbold=cmbx10
  \font\ninebold=cmbx9
  \font\sevenbold=cmbx7
  \font\fivebold=cmbx5

  \input amssym.def \input amssym
     \font\titlemsa=msam10 at 14.4pt
     \font\titlemsb=msbm10 at 14.4pt
     \font\titleeufm=eufm10 at 14.4pt
     \font\twelvemsa=msam10 scaled 1200
     \font\twelvemsb=msbm10 scaled 1200
     \font\twelveeufm=eufm10 scaled 1200
     \font\ninemsa=msam9
     \font\ninemsb=msbm9
     \font\nineeufm=eufm9

   \ifx\cyrfam\undefined
   \else
     \immediate\write16{}%
     \message{ !!! cyr fonts already defined. !!! }
     \message{ --- edit out superfluous font defs? }
   \fi
   \newfam\cyrfam
       \font\titlecyr=wncyr10 scaled 1440 
       \font\twelvecyr=wncyr10 scaled 1200
       \font\tencyr=wncyr10
       \font\ninecyr=wncyr9
       \font\sevencyr=wncyr7
       \font\sixcyr=wncyr6

   \newfam\eusmfam
       \font\titleeusm=eusm10 scaled 1440
       \font\twelveeusm=eusm10 scaled 1200
       \font\teneusm=eusm10
       \font\nineeusm=eusm9
       \font\seveneusm=eusm7
       
       \font\fiveeusm=eusm5

\let\Cal\cal

    \font\ninemrm=cmr9 
    \font\ninei=cmmi9
    \font\ninesy=cmsy9 
    \skewchar\ninei='177
    \skewchar\ninesy='60

  \font\twelvemrm=cmr10 at 12pt 
  \font\twelvei=cmmi10 at 12pt
  \font\twelvesy=cmsy10 at 12pt

  \font\titlemrm=cmr10 at 14.4pt 
  \font\titlei=cmmi10 at 14.4pt
  \font\titlesy=cmsy10 at 14.4pt


  \def\bi{\tenbi}

  \def\Smallfonts{\ninepoint}

  \def\Hfont{\titlepoint\bf}
  \def\Authorfont{\twelvepoint\it}
  \def\HHfont{\twelvepoint\bf}
  \def\HHHfont{\bf}
  \def\Bibfont{\tenbf}
  \def\Coordfont{\nineit }

  \def \thfont {\bf }
  \def \pffont {\it\itSpacing }
  \def \rkfont {\bf }
  \def \dffont {\bf }
  \def \egfont {\bf }

 \def\ninepoint{%
  \def\rm{\fam0\ninerm}%
    \textfont0=\ninemrm  \scriptfont0=\sevenrm  \scriptscriptfont0=\fiverm
    \textfont1=\ninei    \scriptfont1=\seveni   \scriptscriptfont1=\fivei
  \def\mit{\fam1\ninei}%
  \def\oldstyle{\fam1\ninei}%
    \textfont2=\ninesy   \scriptfont2=\sevensy  \scriptscriptfont2=\fivesy
    \textfont3=\tenex    \scriptfont3=\tenex    \scriptscriptfont3=\tenex
  \def\it{\fam\itfam\nineit}%
    \textfont\itfam=\nineit
  \def\bf{\ifmmode\fam\bffam\else\ninebf\fi}%
    \textfont\bffam=\ninebold 
    \scriptfont\bffam=\sevenbold 
    \scriptscriptfont\bffam=\fivebold%
  \def\msa{\fam\msafam\ninemsa}%
    \textfont\msafam=\ninemsa 
    \scriptfont\msafam=\sevenmsa
    \scriptscriptfont\msafam=\fivemsa%
  \def\msb{\fam\msbfam\ninemsb}%
    \textfont\msbfam=\ninemsb%
    \scriptfont\msbfam=\sevenmsb%
    \scriptscriptfont\msbfam=\fivemsb%
  \def\eufm{\fam\eufmfam\nineeufm}%
    \textfont\eufmfam=\nineeufm
    \scriptfont\eufmfam=\seveneufm
    \scriptscriptfont\eufmfam=\fiveeufm
   \def\eusm{\fam\eusmfam\nineeusm}%
     \textfont\eusmfam=\nineeusm
     \scriptfont\eusmfam=\seveneusm
     \scriptscriptfont\eusmfam=\fiveeusm
   \def\cyr{\fam\cyrfam\ninecyr}%
     \textfont\cyrfam=\ninecyr
     \scriptfont\cyrfam=\sevencyr
     \scriptscriptfont\cyrfam=\sixcyr
  \setbox\strutbox=\hbox{\vrule
      height7pt depth3pt width0pt}%
   \baselineskip=10.8pt\rm}

 \let\eightpoint\ninepoint 

 \def\tenpoint{%
  \def\rm{\fam0\tenrm}%
    \textfont0=\tenmrm \scriptfont0=\sevenrm \scriptscriptfont0=\fiverm%
  \def\mit{\fam1\teni}%
  \def\oldstyle{\fam1\teni}%
    \textfont1=\teni   \scriptfont1=\seveni  \scriptscriptfont1=\fivei%
    \textfont2=\tensy  \scriptfont2=\sevensy \scriptscriptfont2=\fivesy%
    \textfont3=\tenex  \scriptfont3=\tenex   \scriptscriptfont3=\tenex%
  \def\it{\fam\itfam\tenit}%
    \textfont\itfam=\tenit%
  \def\bf{\ifmmode\fam\bffam\else\tenbf\fi}%
    \textfont\bffam=\tenbold
    \scriptfont\bffam=\sevenbold%
    \scriptscriptfont\bffam=\fivebold%
  \def\msa{\fam\msafam\tenmsa}%
    \textfont\msafam=\tenmsa%
    \scriptfont\msafam=\sevenmsa%
    \scriptscriptfont\msafam=\fivemsa%
  \def\msb{\fam\msbfam\tenmsb}%
    \textfont\msbfam=\tenmsb%
    \scriptfont\msbfam=\sevenmsb%
    \scriptscriptfont\msbfam=\fivemsb%
  \def\eufm{\fam\eufmfam\teneufm}%
   \textfont\eufmfam=\teneufm
   \scriptfont\eufmfam=\seveneufm
   \scriptscriptfont\eufmfam=\fiveeufm
   \def\eusm{\fam\eusmfam\teneusm}%
    \textfont\eusmfam=\teneusm
    \scriptfont\eusmfam=\seveneusm
    \scriptscriptfont\eusmfam=\fiveeusm
   \def\cyr{\fam\cyrfam\tencyr}%
    \textfont\cyrfam=\tencyr
    \scriptfont\cyrfam=\sevencyr
    \scriptscriptfont\cyrfam=\sixcyr
  \setbox\strutbox=\hbox{\vrule %
      height8.5pt depth3.5ptwidth0pt}%
  \baselineskip=\StdBaselineskip\rm}

 \def\twelvepoint{%
  \def\rm{\fam0\twelverm}%
    \textfont0=\twelvemrm \scriptfont0=\tenmrm \scriptscriptfont0=\sevenrm
    \textfont1=\twelvei   \scriptfont1=\teni   \scriptscriptfont1=\seveni
  \def\mit{\fam1\twelvei}%
  \def\oldstyle{\fam1\twelvei}%
    \textfont2=\twelvesy  \scriptfont2=\tensy  \scriptscriptfont2=\sevensy
    \textfont3=\tenex  \scriptfont3=\tenex  \scriptscriptfont3=\tenex
  \def\it{\fam\itfam\twelveit}%
    \textfont\itfam=\twelveit
  \def\bf{\ifmmode\fam\bffam\else\twelvebf\fi}%
    \textfont\bffam=\twelvebold
    \scriptfont\bffam=\tenbold%
    \scriptscriptfont\bffam=\sevenbold%
  \def\msa{\fam\msafam\twelvemsa}%
    \textfont\msafam=\twelvemsa%
    \scriptfont\msafam=\tenmsa%
    \scriptscriptfont\msafam=\sevenmsa%
  \def\msb{\fam\msbfam\twelvemsb}%
    \textfont\msbfam=\twelvemsb%
    \scriptfont\msbfam=\tenmsb%
    \scriptscriptfont\msbfam=\sevenmsb%
  \def\eufm{\fam\eufmfam\twelveeufm}%
   \textfont\eufmfam=\twelveeufm
   \scriptfont\eufmfam=\teneufm
   \scriptscriptfont\eufmfam=\seveneufm
   \def\eusm{\fam\eusmfam\twelveeusm}%
    \textfont\eusmfam=\twelveeusm
    \scriptfont\eusmfam=\teneusm
    \scriptscriptfont\eusmfam=\seveneusm
   \def\cyr{\fam\cyrfam\tencyr}%
    \textfont\cyrfam=\twelvecyr
    \scriptfont\cyrfam=\tencyr
    \scriptscriptfont\cyrfam=\sevencyr
  \setbox\strutbox=\hbox{\vrule
      height10.2pt depth4.55pt width0pt}%
  \baselineskip=14pt\rm}

 \def\titlepoint{%
    \textfont0=\titlemrm \scriptfont0=\twelvemrm \scriptscriptfont0=\tenmrm
    \textfont1=\titlei   \scriptfont1=\twelvei   \scriptscriptfont1=\teni
  \def\mit{\fam1\titlei}%
  \def\oldstyle{\fam1\titlei}%
    \textfont2=\titlesy  \scriptfont2=\twelvesy  \scriptscriptfont2=\tensy
    \textfont3=\tenex
    \scriptfont3=\tenex
    \scriptscriptfont3=\tenex
  \def\it{\fam\itfam\titleit}%
    \textfont\itfam=\titleit
  \def\bf{\ifmmode\fam\bffam\else\titlebf\fi}%
    \textfont\bffam=\titlebold
    \scriptfont\bffam=\twelvebold%
    \scriptscriptfont\bffam=\tenbold%
  \def\msa{\fam\msafam\titlemsa}%
    \textfont\msafam=\titlemsa%
    \scriptfont\msafam=\twelvemsa%
    \scriptscriptfont\msafam=\tenmsa%
  \def\msb{\fam\msbfam\titlemsb}%
    \textfont\msbfam=\titlemsb%
    \scriptfont\msbfam=\twelvemsb%
    \scriptscriptfont\msbfam=\tenmsb%
  \def\eufm{\fam\eufmfam\titleeufm}%
    \textfont\eufmfam=\titleeufm
    \scriptfont\eufmfam=\twelveeufm
    \scriptscriptfont\eufmfam=\teneufm
   \def\eusm{\fam\eusmfam\titleeusm}%
     \textfont\eusmfam=\titleeusm
     \scriptfont\eusmfam=\twelveeusm
     \scriptscriptfont\eusmfam=\teneusm
   \def\cyr{\fam\cyrfam\tencyr}%
    \textfont\cyrfam=\titlecyr
    \scriptfont\cyrfam=\twelvecyr
    \scriptscriptfont\cyrfam=\tencyr
  \setbox\strutbox=\hbox{\vrule
      height12.3pt depth5.54pt width0pt}%
  \baselineskip=16pt\rm}

\newbox\AuthorBox\newbox\TitleBox
\newbox\TFLinebox
\newbox\FLinebox
\newbox\HLinebox
\def\SetTFLinebox#1{\setbox\TFLinebox=\hbox{#1}}
\def\SetFLinebox#1{\setbox\FLinebox=\hbox{#1}}
\def\SetHLinebox#1{\setbox\HLinebox=\hbox{#1}}

 \def\SetAuthorHead#1{%
     \setbox\AuthorBox=\hbox{\ninepoint \it 
           \ignorespaces\frenchspacing#1\unskip}}
 \def\SetTitleHead#1{%
     \setbox\TitleBox=\hbox{\ninepoint \it
           \ignorespaces\frenchspacing#1\unskip}}

  \def\itSpacing{\relax}
  \def\itSpacingOff{\relax}


 \def\Hrule{\hrule width0pt height0pt}

  \newskip\ProcSkip \ProcSkip 8pt plus2pt minus2pt

 \newskip\LastSkip
 \def\SaveLastSkip{\LastSkip\lastskip}
 \def\RestoreLastSkip{\vskip-\LastSkip\vskip\LastSkip}

 \def\NoindentAfter{\everypar={\setbox0=\lastbox\everypar={}}}

 \long\def\H#1\par#2\par{\notenumber=0 \titlepagetrue%
    {
    \baselineskip=20pt
    \parindent=0pt\parskip=0pt\frenchspacing
    \leftskip=0pt plus .2\hsize minus .3\hsize
    \rightskip=0pt plus .2\hsize minus .3\hsize
 \def\\{\unskip\break}%
    \pretolerance=10000 \Hfont #1\unskip\break
     \vskip7pt\Hrule
\hfill \Authorfont #2\hfill\hfill\unskip}
    \vskip48pt plus 4pt minus 4pt
    \par\NoindentAfter\rm}

 \long\def\Hi#1\par#2\par{\notenumber=0 \titlepagetrue%
    {  \baselineskip=0pt  \parindent=0pt\parskip=0pt\frenchspacing
    \leftskip=0pt plus .2\hsize minus .3\hsize
    \rightskip=0pt plus .2\hsize minus .3\hsize
}
    \rm}


 \newdimen\PageRemainder
  \def\SetPageRemainder{
     \PageRemainder=\pagegoal
     \ifdim\PageRemainder=\maxdimen\PageRemainder=\vsize
     \else\advance\PageRemainder by -1\pagetotal\fi}

  \def\Rpt@{}\def\Rpt@@{}

  \long\def\HH#1\par{\par
  \SaveLastSkip\removelastskip\goodbreak
  \ifdim\LastSkip<30pt 
     \LastSkip 30pt
plus 3pt minus 2pt\fi
  \SetPageRemainder\advance\PageRemainder-\LastSkip
  \ifdim\PageRemainder<150pt
       \edef\Rpt@{remain = \the\PageRemainder\noexpand\\
                pagetotal=\the\pagetotal\noexpand\\
                           pagegoal=\the\pagegoal}%
          \fi
   \ifdim\PageRemainder<65pt 
       \ifdim\PageRemainder > 0pt
          \edef\Rpt@@{\noexpand\\
                      Had HH PageRemainder$<$\relax 65pt\noexpand\\
                      Hence forced break!}%
     \vskip 0pt plus .2\PageRemainder\eject 
    \fi\fi
    \vskip\LastSkip\Hrule 
    \pretolerance=10000\rightskip=0pt plus 3em
    \hangafter1 \hangindent=2.2em%
    \noindent
    \HHfont \unskip \Ednote{\Rpt@\Rpt@@}%
            \def\Rpt@{}\def\Rpt@@{}%
            \ignorespaces
            #1\par\rightskip=0pt\pretolerance=\StdPretolerance%
    \NoindentAfter
\tenpoint\rm%
     \medskip \vskip\ProcSkip}

  \long\def\HHH#1\par{\par%
  \SaveLastSkip\removelastskip\goodbreak
  \ifdim\LastSkip<\ProcSkip%
     \LastSkip\ProcSkip\fi
  \SetPageRemainder\advance\PageRemainder-\LastSkip
  \ifdim\PageRemainder<150pt
       \edef\Rpt@{remain = \the\PageRemainder\noexpand\\
                pagetotal=\the\pagetotal\noexpand\\
                           pagegoal=\the\pagegoal}%
       \fi
   \ifdim\PageRemainder<48pt  
        \ifdim\PageRemainder > 0pt
             \edef\Rpt@@{\noexpand\\
                      Had HHH PageRemainder$<$\relax48pt\noexpand\\
                      Hence forced break!}%
       \vskip 0pt plus .2\PageRemainder\eject 
      \fi\fi
   \vskip\LastSkip\par\noindent
   \HHHfont \unskip\Ednote{\Rpt@\Rpt@@}%
  \def\Rpt@{}\def\Rpt@@{}%
  \ignorespaces
   #1\unskip.\quad\rm\ignorespaces
   \ignorepars}

  \long\def\ignorepars#1\par{\def\Test{#1}%
     \ifx\Test\Empty\def\This{\ignorepars}%
        \else\def\This{\Test\par}\fi
           \This}
  \def\Empty{}

 \def\Abstract#1\par{\bgroup\Smallfonts\narrower\HHH #1\par}
 \def\endAbstract{\par\egroup}


 \def\ProcBreak{\par%
    \ifdim\lastskip<8pt%
    \removelastskip%
    \penalty-200\vskip\ProcSkip\fi}

 \def\th#1\par{\ProcBreak \noindent
   {\thfont\ignorespaces
    #1\unskip.}\it\itSpacing\kern.4em\ignorepars}

 \def\endth{\ProcBreak\rm\itSpacingOff }


 \def\pf#1\par{\ProcBreak %
    \noindent\pffont#1\unskip.\rm\itSpacingOff{\kern .7em}\ignorepars}

 \def\endpf{\medskip \ProcBreak } 

  \def\qedbox{\hbox{\vbox{
    \hrule width0.2cm height0.2pt
    \hbox to 0.2cm{\vrule height 0.2cm width 0.2pt
             \hfil\vrule height0.2cm width 0.2pt}
    \hrule width0.2cm height 0.2pt}\kern1pt}}

  \def\qed{\ifmmode\qedbox
    \else\unskip\ \hglue0mm\hfill\qedbox\ProcBreak\fi}

  \def \rk #1\par{\ProcBreak
     \noindent{\rkfont\ignorespaces #1\unskip.}%
     \rm\kern.6em\ignorepars}

  \def \endrk {\medskip\ProcBreak }

  \def \df #1\par{\ProcBreak
     \noindent{\dffont\unskip\ignorespaces #1\unskip.}%
     \rm\kern.6em\ignorepars}

  \def \enddf {\medskip\ProcBreak }

  \def \eg #1\par{\ProcBreak
     \noindent\egfont\unskip\ignorespaces #1\unskip.
     \rm\kern.6em\ignorepars}

  \def \endeg {\medskip\ProcBreak }

  \newdimen\Overhang

   \def\MaxTag@#1#2#3#4#5{\setbox0=\hbox{#4\ignorespaces#2\unskip}%
     \dimen0=\wd0\advance\dimen0 by#3
     \ifdim\dimen0<#5\relax\dimen0=#5\fi
     \expandafter\edef\csname #1Hang\endcsname{\the\dimen0}}

 \def\MaxItemTag#1{\MaxTag@{Item}{#1}{.4em}{\ItemStyle}{\parindent}}%
 \def\MaxItemItemTag#1{%
        \MaxTag@{ItemItem}{#1}{.4em}{\ItemItemStyle}{\parindent}}
 \def\MaxNrTag#1{\MaxTag@{Nr}{#1}{.5em}{\NrStyle}{\parindent}}
 \def\MaxReferenceTag#1{%
        \MaxTag@{Reference}{[#1]}{.6em}{\ninerm}{\parindent}}
 \def\MaxFootTag#1{\MaxTag@{Foot}{#1}{.4em}{\ninerm}{\z@}}

  \def\SetOverhang@{\Overhang=.8\dimen0%
     \advance\Overhang by \wd0\relax
     \ifdim\Overhang>\hangindent\relax
       \advance\Overhang by .25\dimen0%
       \Ednote{Tag is pushing text.}\osumess{Tag is pushing text.}%
     \else\Overhang=\hangindent
     \fi}

   \def\Item#1{\par\noindent
      \hangafter1\hangindent=\ItemHang
      \setbox0=\hbox{\ItemStyle\ignorespaces#1\unskip}%
      \dimen0=.4em\SetOverhang@
      \rlap{\box0}\kern\Overhang\ignorespaces}

   \def\ItemItem#1{\par\noindent
      \hangafter1\hangindent=\ItemItemHang
      \setbox0=\hbox{\ItemItemStyle\ignorespaces#1\unskip}%
      \dimen0=.4em\SetOverhang@
      \advance\hangindent by \ItemHang
      \kern\ItemHang\rlap{\box0}%
      \kern\Overhang\ignorespaces}

  \def\Nr#1{\par\noindent\hangindent=\NrHang 
    \setbox0=\hbox{\NrStyle\ignorespaces#1\unskip}%
    \dimen0=.5em\SetOverhang@
    \rlap{\box0}\kern\Overhang
    \hangindent=\z@\ignorespaces}

   \newskip\Rosterskip\Rosterskip 1pt plus1pt 
   \def\Roster{\par\ifdim\lastskip<\Rosterskip\removelastskip\vskip\Rosterskip\fi
    \bgroup}
   \def\endRoster{\par\global\edef\LastSkip@{\the\lastskip}\removelastskip
       \egroup\penalty-50\LastSkip\LastSkip@\relax
       \ifdim\LastSkip<\Rosterskip\LastSkip\Rosterskip\fi
       \vskip\LastSkip}


 \long\def \Em #1{{%
   \bi#1\unskip\/}}



 \def\cite#1{
    \def\nextiii@##1,##2\end@{{\frenchspacing\rm 
      \lBr\ignorespaces##1\unskip{\rm,~\ignorespaces##2}\rBr}}%
    \IN@0,@#1@%
    \ifIN@\def\next{\nextiii@#1\end@}\else
    \def\next{{\rm\lBr#1\rBr}}\fi\next}


   \def \Bib#1\par{%
       \par\removelastskip\SetPageRemainder
       \ifdim\PageRemainder < 97pt
        \ifdim\PageRemainder > 0pt
        \vfill\eject
       \fi\fi
    \ProcBreak \par\begingroup\parskip=0 pt%
    \goodbreak \vskip 15 pt plus 10 pt
    \noindent\null\hfill\Bibfont
      \ignorespaces #1\unskip\hfill\null\par 
    \frenchspacing \Smallfonts\rm
    \parskip=2.5 pt plus 1 pt minus.5pt%
    \nobreak\vskip 12pt plus 2pt minus2pt\nobreak
    \leftskip=0 pt \baselineskip=10.5pt}

 \def\ReferenceTagSlide{0em}
  \def\ReferenceTagGap{.5em}

  \def \rf#1{\par\noindent
     \hangafter1\hangindent=\ReferenceHang      
     \setbox0=\hbox{\ninerm[\ignorespaces#1\unskip]}%
     \dimen0=\ReferenceTagGap\SetOverhang@
     \rlap{\kern\ReferenceTagSlide\box0}%
     \kern\Overhang\ignorespaces}

  \def\ref#1\par#2\par#3\par#4\par{%
     \rf{#1}#2\unskip,\ #3\unskip,\
     #4\unskip.}

  \def\endBib{\par\endgroup\vskip 12pt minus 6pt }


  \long\def\Coordinates#1\endCoordinates{
 {\par\vskip4pt\def\\{\unskip, }\Coordfont\baselineskip10.5pt\noindent#1}}

 \def\pagecontents{
  \gdef\Pagetot@l{\pagetotal}
  \ifvoid\TRMargIns\else
    \rlap{\kern\hsize\kern10pt\vbox to 0pt{%
         \box\TRMargIns\vss}}\fi
  \ifvoid\topins\else\unvbox\topins\fi
   \dimen@=\dp\@cclv \unvbox\@cclv 
   \ifvoid\footins\else 
     \vskip\skip\footins
     \footnoterule
     \unvbox\footins\fi
   \ifr@ggedbottom \kern-\dimen@ \vfil \fi}


 \newcount\Ht 

 \def \Acc{\expandafter } 

 \def\swthat{\raise -1.1 ex\hbox{\sam$\widehat{}$}}
 \def\swttilde{\raise -1.2 ex\hbox{\sam$\widetilde{}$}}
 \def \overdot{{\raise .2 ex \hbox to 0pt {\hss\bf\smash{.}\hss}}}
 \def \overcircle{{\raise .1 ex \hbox to 0pt
    {\sam$\eightpoint\scriptstyle\hss\circ\hss$}}}

 \def \Mathaccent#1#2{{\sam 
  \setbox4=\hbox{$\vphantom{#2}$}
  \Ht=\ht4 
  \setbox5=\hbox{${#1}$}
  \setbox6=\hbox{${#2}$}
  \setbox7=\hbox to .5\wd6{}
  \copy7\kern .1\Ht \raise\Ht sp\hbox{\copy5}\kern-.1\Ht
  \copy7\llap{\box6}
  }}

  \def\SwtCheck #1{
        \ifmmode \check{#1}%
                \else \v {#1}%
                \fi}

 \def\barpartial {%
   \kern .17 em
    \overline {\kern -.17 em\partial\kern-.03 em}%
    \kern .03 em}

 
  \def\Overline#1{\setbox1=\hbox{\sam ${#1}$}%
      \ifdim \wd1 > 6pt
    \kern .11 em
    \overline {\kern -.11 em#1\kern-.14 em}
    \kern .14 em
  \else
    \kern .03 em
    \overline {\kern -.03 em#1\kern-.04 em}
    \kern .04 em
  \fi}

 \def\SOverline#1{\setbox1=\hbox{\sam ${#1}$}%
      \ifdim \wd1 > 7pt
    \kern .22 em
    \overline {\kern -.22 em#1\kern-.09 em}%
    \kern .09 em
  \else
    \kern .10 em
    \overline {\kern -.10 em#1\kern-.04 em}%
    \kern .04 em
  \fi}


 \def\Underline#1{\setbox1=\hbox{\sam ${#1}$}%
      \ifdim \wd1 > 6pt
    \kern .11 em
    \underline {\kern -.11 em#1\kern-.14 em}
    \kern .14 em
  \else
    \kern .03 em
    \underline {\kern -.03 em#1\kern-.04 em}
    \kern .04 em
  \fi}

 \def\SUnderline#1{\setbox1=\hbox{\sam ${#1}$}%
      \ifdim \wd1 > 7pt
    \kern .04 em
    \underline {\kern -.04 em#1\kern-.2 em}%
    \kern .2 em
  \else
    \kern .0 em
    \underline {\kern -.0 em#1\kern-.15 em}%
    \kern .15 em
  \fi}


 \def \Blackbox
   {\leavevmode\hskip .3pt \vbox
   {\hrule height 5pt\hbox{\hskip 4.5pt}}\hskip .5pt}

 \def \XX{\Blackbox\kern.5pt\Blackbox} 

  \def\.{.\kern1pt}

    \def\Hyphen{\edef\this{\the\hyphenchar\font}%
          \hyphenchar\font=-1\char\this\hyphenchar\font=\this}

 \ifx\undefined\text
  \def\text#1{\hbox{\rm #1}}\fi 



   \everymath{}  

  \def\PassMath@@{\aftergroup\AfterMath@} 

 \let\PassMath@\PassMath@@

 \def\AfterMath@{\futurelet\next\AfterMathMole@}

 \def\AfterMathMole@{
      \ifcat\next\space
          \def\this{}
      \else
      \ifcat\next\egroup %
        \def\this{\osumess{Handset mathsurround?? ---(see dollar brace)}}%
      \else
      \def\this{\AAfterMath@}
      \fi\fi
      \this}

 \def\hyphen@{-}
 \def\paren@{)}
 \def\apostr@{'}

 \def\MSC#1{\kern-.8\mathsurround#1\kern.8\mathsurround}

 \def\AAfterMath@#1{\def\Next{#1}
    \IN@0\Next @,.;:!?\relax @%
    \ifIN@\def\this{\MSC{\Next}}%
    \else
    \ifx\Next\hyphen@\def\this{\futurelet\next\AfterHyphen@}%
    \else
    \ifx\Next\paren@\def\this{#1}%
    \else 
    \ifx\Next\apostr@\def\this{#1}%
    \else \def\this{\osumess{Handset mathsurround??}%
                 #1}\fi\fi\fi\fi
    \this}

 \def\AfterHyphen@#1{\def\Next{#1}%
   \ifx\Next\hyphen@\def\this{--}\else
   \ifcat\next\space%
   \def\this{\kern-\mathsurround\kern.05em- \Next}\else
   \def\this{\kern-\mathsurround\kern.05em\Hyphen\Next}\fi\fi\this}

 \def\sam{\mathsurround=\z@\let\PassMath@\relax}  %
 \def\mas{\mathsurround=\StdMathsurround\let\PassMath@\PassMath@@}
 
 \def\Mas{\mathsurround=\StdMathsurround
                \everymath{\PassMath@}\let\PassMath@\PassMath@@}

 \def\m@th{\mathsurround=\z@\everymath{}}

 \def\m@@th{\mathsurround=\z@\everymath={}\let\m@th\relax}

\def\underbar#1{$\setbox\z@\hbox{#1}\dp\z@\z@
      \m@th \underline{\box\z@}$\relax}

\def\mathhexbox#1#2#3{\leavevmode
  \hbox{\m@@th$\m@th \mathchar"#1#2#3$}}

\def\dots{\relax\ifmmode\ldots\else$\m@th\ldots\,$\relax\fi}

\def\dotfill{\cleaders\hbox{\m@@th$\m@th \mkern1.5mu.\mkern1.5mu$}\hfill}
\def\rightarrowfill{$\m@th\mathord-\mkern-6mu%
  \cleaders\hbox{\m@@th$\mkern-2mu\mathord-\mkern-2mu$}\hfill
  \mkern-6mu\mathord\rightarrow$\relax}
\def\leftarrowfill{$\m@th\mathord\leftarrow\mkern-6mu%
  \cleaders\hbox{\m@@th$\mkern-2mu\mathord-\mkern-2mu$}\hfill
  \mkern-6mu\mathord-$\relax}

\def\downbracefill{$\m@th\braceld\leaders\vrule\hfill\braceru
  \bracelu\leaders\vrule\hfill\bracerd$\relax}
\def\upbracefill{$\m@th\bracelu\leaders\vrule\hfill\bracerd
  \braceld\leaders\vrule\hfill\braceru$\relax}

\def\angle{{\vbox{\m@@th\ialign{$\m@th\scriptstyle##$\crcr
      \not\mathrel{\mkern14mu}\crcr
      \noalign{\nointerlineskip}
      \mkern2.5mu\leaders\hrule height.34pt\hfill\mkern2.5mu\crcr}}}}

\def\big#1{{\m@@th\hbox{$\left#1\vbox to8.5\p@{}\right.\n@space$}}}
\def\Big#1{{\m@@th\hbox{$\left#1\vbox to11.5\p@{}\right.\n@space$}}}
\def\bigg#1{{\m@@th\hbox{$\left#1\vbox to14.5\p@{}\right.\n@space$}}}
\def\Bigg#1{{\m@@th\hbox{$\left#1\vbox to17.5\p@{}\right.\n@space$}}}
\def\n@space{\nulldelimiterspace\z@ \m@th}

\def\root#1\of{\setbox\rootbox\hbox{\m@@th$\m@th\scriptscriptstyle{#1}$}
  \mathpalette\r@@t}
\def\r@@t#1#2{\setbox\z@\hbox{\m@@th$\m@th#1\sqrt{#2}$\relax}
  \dimen@\ht\z@ \advance\dimen@-\dp\z@
  \mkern5mu\raise.6\dimen@\copy\rootbox \mkern-10mu \box\z@}

\def\mathph@nt#1#2{\setbox\z@\hbox{\m@@th$\m@th#1{#2}$}\finph@nt}

\def\mathsm@sh#1#2{\setbox\z@\hbox{\m@@th$\m@th#1{#2}$}\finsm@sh}

\def\@vereq#1#2{\lower.5\p@\vbox{\m@@th\baselineskip\z@skip\lineskip-.5\p@
    \ialign{$\m@th#1\hfil##\hfil$\crcr#2\crcr=\crcr}}}

\def\mathpalette#1#2{\sam\mathchoice{#1\displaystyle{#2}}%
  {#1\textstyle{#2}}{#1\scriptstyle{#2}}{#1\scriptscriptstyle{#2}}\mas}

\def\widehat#1{\setbox\z@\hbox{\sam$#1$}%
 \ifdim\wd\z@>\tw@ em\mathaccent"0\msbfam@5B{#1}%
 \else\mathaccent"0362{#1}\fi}
\def\widetilde#1{\setbox\z@\hbox{\sam$#1$}%
 \ifdim\wd\z@>\tw@ em\mathaccent"0\msbfam@5D{#1}%
 \else\mathaccent"0365{#1}\fi}

 \def\dots{\relax{}
  \ifmmode\def\thedots{\mdots@}\else\def\thedots{\tdots@}\fi %
  \thedots}

 \let\@ldeqno\eqno\let\@ldleqno\leqno
 \def\eqno{\everymath{}\@ldeqno} \def\leqno{\everymath{}\@ldleqno}

  \let\@ldeqalignno\eqalignno
  \def\eqalignno#1{\sam\@ldeqalignno{#1}\mas}
  \let\@ldeqalign\eqalign
  \def\eqalign#1{\sam\@ldeqalign{#1}\mas}

 \def\overrightarrow#1{\vbox{\m@th\ialign{##\crcr
      \rightarrowfill\crcr\noalign{\kern-\p@\nointerlineskip}
      $\hfil\displaystyle{#1}\hfil$\crcr}}}
 \def\overleftarrow#1{\vbox{\m@th\ialign{##\crcr
      \leftarrowfill\crcr\noalign{\kern-\p@\nointerlineskip}
      $\hfil\displaystyle{#1}\hfil$\crcr}}}
 \def\overbrace#1{\mathop{\vbox{\m@th\ialign{##\crcr\noalign{\kern3\p@}
      \downbracefill\crcr\noalign{\kern3\p@\nointerlineskip}
      $\hfil\displaystyle{#1}\hfil$\crcr}}}\limits}
 \def\underbrace#1{\mathop{\vtop{\m@th\ialign{##\crcr
      $\hfil\displaystyle{#1}\hfil$\crcr\noalign{\kern3\p@\nointerlineskip}
      \upbracefill\crcr\noalign{\kern3\p@}}}}\limits}

  \let\@ldmatrix\matrix
  \let\end@ldmatrix\endmatrix
  \def\matrix{\sam\@ldmatrix}
  \def\endmatrix{\end@ldmatrix\mas}
  \let\@ldgather\gather
  \let\end@ldgather\endgather
  \def\gather{\sam\@ldgather}
  \def\endgather{\end@ldgather\mas}
  \let\@ldalign\align
  \let\end@ldalign\endalign
  \def\align{\sam\@ldalign}
  \def\endalign{\end@ldalign\mas}
  \let\@ldaligned\aligned
  \let\end@ldaligned\endaligned
  \def\aligned{\sam\@ldaligned}
  \def\endaligned{\end@ldaligned\mas}
  \let\@ldtag\tag
  \def\tag{\sam\@ldtag}
   %

   \let\MinCDArrowWidth\minCDaw@




\newskip\insertskipamount\newskip\inserthardskipamount
\insertskipamount 6pt plus2pt 
\inserthardskipamount 6pt
\def\insertskip{\vskip\insertskipamount}
\newcount\SplitTest
\def\SetSplitTest{\SplitTest\insertpenalties
  \insert\topins{\floatingpenalty1}%
  \advance\SplitTest-\insertpenalties}
\def\midinsert{\par
 \SaveLastSkip\penalty-150\SetSplitTest\RestoreLastSkip
 \ifnum\SplitTest=-1
  \@midfalse\p@gefalse\else\@midtrue\fi\@ins}
\def\@ins{\par\begingroup\setbox\z@\vbox\bgroup%
  \vglue\inserthardskipamount}
\def\endinsert{\egroup 
  \if@mid \dimen@\ht\z@ \advance\dimen@\dp\z@
    \advance\dimen@\insertskipamount
    \advance\dimen@\pagetotal\advance\dimen@-\pageshrink
    \ifdim\dimen@>\pagegoal\@midfalse\p@gefalse\fi\fi
  \if@mid%
    \ifdim\lastskip<\insertskipamount\removelastskip\insertskip\fi
    \nointerlineskip\box\z@\penalty-200\insertskip
  \else%
    \SaveLastSkip
    \insert\topins{\penalty100 
    \splittopskip\z@skip
    \splitmaxdepth\maxdimen \floatingpenalty\z@
    \ifp@ge \dimen@\dp\z@
    \vbox to\vsize{\unvbox\z@\kern-\dimen@}
    \else \box\z@\nobreak\insertskip\fi}
    \RestoreLastSkip
   \fi\endgroup}


  \newcount\notenumber
  
  \def\note{\advance\notenumber by 1
    \footnote{\the\notenumber)}}

  \newbox\footbox

  \def\footnote#1{\let\@sf\empty
    \ifhmode\edef\@sf{\spacefactor\the\spacefactor}\/\fi
    \sam${}^{\fam0 #1}$\@sf\vfootnote{#1}}%

  \def\vfootnote#1{\insert\footins\bgroup
     \interlinepenalty100 \splittopskip=1pt
     \floatingpenalty=20000
     \leftskip=0pt\rightskip=0pt%
     \parindent=.3em
     \Smallfonts\rm
     \FootItem@{#1}
     \futurelet\next\fo@t}

  \def\FootItem@#1{\par\hangafter1\hangindent=\FootHang
     \setbox0=\hbox{\ignorespaces#1\unskip}%
     \dimen0=.4em\SetOverhang@
     \noindent\rlap{\box0}\kern\Overhang\ignorespaces}


  \def\fo@t{\ifcat\bgroup\noexpand\next \let\next\f@@t
    \else\let\next\f@t\fi \next}
  \def\f@@t{\bgroup\aftergroup\@foot\let\next}
  \def\f@t#1{\baselineskip=10pt\lineskip=1pt
            \lineskiplimit=0pt #1\@foot}%
  \def\@foot{
        \hbox{\vrule height0pt depth5pt width0pt}
        \egroup}
  \skip\footins=12 pt plus 0pt minus 0pt 
  \count\footins=1000 
  \dimen\footins=8in 



 \def\osumess#1{\EdSpider{\immediate\write16{Line \the\inputlineno: #1}}}%
 \def\HideEdStuff{\gdef\EdSpider##1{}}

 \font\BigSym=cmmi10 scaled \magstep 4

 \def\change{\InLMargin{\hbox{\BigSym \char63\kern10pt}}}

 \def\beginchange{\InLMargin{\hbox{\sam\twelvepoint$\heartsuit$\kern10pt}}}

 \def\endchange{\InLMargin{\hbox{\sam\twelvepoint$\spadesuit$\kern10pt}}}

 \def\InLMargin#1{\strut\vadjust{%
     \kern-\strutdepth
     \vtop to \strutdepth{%
         \baselineskip\strutdepth
         \llap{\sam$\smash{\hbox{\EdSpider{#1}}}$}\null}}}

 \def\strutdepth{\dp\strutbox}
 \def\strutheight{\ht\strutbox}

 \def\NoteInRMargin#1{\strut\vadjust{%
     \kern-1.001\strutdepth
     \vtop to \strutdepth{%
       \baselineskip\strutdepth
       \vss\rlap{\ninepoint\unskip\hskip\hsize
         \vtop to 0pt{%
           \hsize=16em\hfuzz=\hsize
           \leftskip=10pt%
           \rightskip=0pt plus 10000pt%
           \baselineskip=9.8pt\lineskip=.2pt%
           \let\\\break
           \noindent\EdSpider{#1}\vss}%
                \kern10pt}\hbox{}}
       }}

 \def\ednote#1{\NoteInRMargin{\tentt #1}}

 \def\cbar{\InLMargin{%
      \dimen0=\strutdepth\advance\dimen0 by \lineskip
      \vrule width 3pt
      height \strutheight depth \dimen0 \kern
      3pt}}

 \def\ccbar{\InLMargin{%
      \dimen0=2\strutdepth\advance\dimen0 by 2\lineskip
      \vrule width 3pt
        height 3\strutheight depth \dimen0 \kern
      3pt}}

 \newinsert\TRMargIns
 \dimen\TRMargIns=\maxdimen

  \def\Ednote#1{\insert\TRMargIns{%
       \vbox to 0pt{\hsize=140pt\hfuzz=\hsize
           \leftskip=6pt%
           \rightskip=0pt plus 10000pt%
           \baselineskip=9.8pt\lineskip=.2pt%
           \let\\\break
           \SetPageRemainder
           \vglue540pt\vglue-\PageRemainder
           \noindent\EdSpider{\tentt #1}\vss}%
       \smallskip}}

 \def\KillEdStuff{\def\ednote##1{}\def\Ednote##1{}%
      \let\change\relax\let\beginchange\relax\let\endchange\relax
       \let\cbar\relax\let\ccbar\relax}


  \topskip=12pt
  \newskip\StdBaselineskip 
  \StdBaselineskip 12pt
  \lineskip=1.1pt
  \lineskiplimit=.8pt
  \widowpenalty=10000 
  \clubpenalty=10000  
  \abovedisplayskip=6pt plus 1pt minus 1pt
  \abovedisplayshortskip=3pt plus 1.5pt
  \belowdisplayskip=6pt plus 1pt minus 1pt
  \belowdisplayshortskip=5pt plus 1pt minus 1pt
  \hfuzz=1.5pt   

  \def\StdPretolerance{100}
  \tolerance=\StdPretolerance

  \newdimen\StdMathsurround
  \StdMathsurround=1.5pt 
  \mathsurround=\StdMathsurround
  \Mas                   

   \def\prose{\relax\hbox{\kern.6\StdMathsurround}}
  
  \def\StdParskip{0pt}    
  \parskip=\StdParskip
  \parindent=0.5cm
 

  \def\Times{ptmr  } 
  \def\TimesI{ptmri  } 
  \def\TimesB{ptmb  }
  \def\TimesBI{ptmbi  }
  \def\HelveticaN{phvrrn }

  =\Times at 10bp
  =\TimesB at 10bp
  \font\tenit=\TimesI at 10bp
  =\TimesBI at 10bp

  \font\tenmrm=cmr10  


    =\Times at 9bp 
    \font\nineit=\TimesI at 9bp 
    =\TimesB at 9bp 
    =\TimesBI at 9bp 

    =\HelveticaN at 9bp 


  =\Times at 12bp
  \font\twelveit=\TimesI at 12bp
  =\TimesB at 12bp


  \font\titleit=\TimesI at 14.4bp
  =\TimesB at 14.4bp

 \SetAuthorHead{AuthorHead} 
 \SetTitleHead{TitleHead}  


  \def\lBr{\raise.125ex\hbox{[\kern.1125ex}}
  \def\rBr{\raise.125ex\hbox{\kern.1125ex]}}

 \setbox\footbox=\hbox{\Smallfonts 2)~}



  \bgroup
  \catcode`\@=11 
  \gdef\itSpacing{%
     \xspaceskip=.31em plus.1em minus.05em \sfcode `f=2001
     \itWarning@\let\itWarning@\itWarning@@}
  \gdef\itSpacingOff{%
     \xspaceskip=0pt \sfcode `f=1000
     \let\itWarning@\relax}
   \global\let\itWarning@\relax
  \gdef\itWarning@@{\errmessage{%
  Special italic spacing already in force
  (you have probably omitted an ``endth'').
  See itSpacing macro in osuPSfnt.sty
         }}
  \egroup

 \fontdimen1\titlebf=0.0pt
 \fontdimen2\titlebf=3.6135pt
 \fontdimen3\titlebf=2.8908pt
 \fontdimen4\titlebf=1.44539pt
 \fontdimen5\titlebf=6.64882pt
 \fontdimen6\titlebf=14.45398pt
 \fontdimen7\titlebf=1.60439pt

 \fontdimen1\tenbi=0.26794pt
 \fontdimen2\tenbi=2.50937pt
 \fontdimen3\tenbi=2.00749pt
 \fontdimen4\tenbi=1.00374pt
 \fontdimen5\tenbi=4.59717pt
 \fontdimen6\tenbi=10.03749pt
 \fontdimen7\tenbi=1.11415pt

 \fontdimen1\twelverm=0.0pt
 \fontdimen2\twelverm=3.01125pt
 \fontdimen3\twelverm=2.409pt
 \fontdimen4\twelverm=1.2045pt
 \fontdimen5\twelverm=5.39615pt
 \fontdimen6\twelverm=12.045pt
 \fontdimen7\twelverm=1.33699pt

 \fontdimen1\twelveit=0.27731pt
 \fontdimen2\twelveit=3.01125pt
 \fontdimen3\twelveit=2.409pt
 \fontdimen4\twelveit=1.2045pt
 \fontdimen5\twelveit=5.37207pt
 \fontdimen6\twelveit=12.045pt
 \fontdimen7\twelveit=1.33699pt

 \fontdimen1\twelvebf=0.0pt
 \fontdimen2\twelvebf=3.01125pt
 \fontdimen3\twelvebf=2.409pt
 \fontdimen4\twelvebf=1.2045pt
 \fontdimen5\twelvebf=5.5407pt
 \fontdimen6\twelvebf=12.045pt
 \fontdimen7\twelvebf=1.33699pt

 \fontdimen1\tenrm=0.0pt
 \fontdimen2\tenrm=2.50937pt
 \fontdimen3\tenrm=2.00749pt
 \fontdimen4\tenrm=1.00374pt
 \fontdimen5\tenrm=4.49678pt
 \fontdimen6\tenrm=10.03749pt
 \fontdimen7\tenrm=1.11415pt

 \fontdimen1\tenit=0.27731pt
 \fontdimen2\tenit=2.50937pt
 \fontdimen3\tenit=2.00749pt
 \fontdimen4\tenit=1.00374pt
 \fontdimen5\tenit=4.47672pt
 \fontdimen6\tenit=10.03749pt
 \fontdimen7\tenit=1.11415pt

 \fontdimen1\tenbf=0.0pt
 \fontdimen2\tenbf=2.50937pt
 \fontdimen3\tenbf=2.00749pt
 \fontdimen4\tenbf=1.00374pt
 \fontdimen5\tenbf=4.61723pt
 \fontdimen6\tenbf=10.03749pt
 \fontdimen7\tenbf=1.11415pt

 \fontdimen1\ninerm=0.0pt
 \fontdimen2\ninerm=2.25842pt
 \fontdimen3\ninerm=1.80673pt
 \fontdimen4\ninerm=0.90337pt
 \fontdimen5\ninerm=4.0471pt
 \fontdimen6\ninerm=9.03374pt
 \fontdimen7\ninerm=1.00273pt

 \fontdimen1\nineit=0.27731pt
 \fontdimen2\nineit=2.25842pt
 \fontdimen3\nineit=1.80673pt
 \fontdimen4\nineit=0.90337pt
 \fontdimen5\nineit=4.02904pt
 \fontdimen6\nineit=9.03374pt
 \fontdimen7\nineit=1.00273pt

 \fontdimen1\ninebf=0.0pt
 \fontdimen2\ninebf=2.25842pt
 \fontdimen3\ninebf=1.80673pt
 \fontdimen4\ninebf=0.90337pt
 \fontdimen5\ninebf=4.15552pt
 \fontdimen6\ninebf=9.03374pt
 \fontdimen7\ninebf=1.00273pt


 \newcount\MaxSpaceFactor
 \MaxSpaceFactor=3000 

 \def\ItemStyle{\rm}
 \def\NrStyle{\rm}
 \def\ItemItemStyle{\rm}

 \MaxItemTag{(iii)}
 \MaxItemItemTag{(iii)}
 \MaxNrTag{(2)}
 \MaxFootTag{2)}
 \def\ReferenceHang{30pt}

 \catcode`\@=\active


\loadbold

=\Times  
=\Times scaled750
=\Times scaled650
\font\rms=\Times scaled 920 

=\TimesBI scaled 860
=\TimesI scaled 860

\textfont0=\rrm  
\scriptfont0=\erm 
\scriptscriptfont0=\srm

\def\Augment#1#2{%
    \toks0\expandafter{#1}\toks2{#2}%
    \edef#1{\the\toks0\the\toks2}}

 \font\twelverma=\Times  scaled 1200
 \font\tenrma=\Times  scaled 1000
 \font\ninerma=\Times scaled 920
 =\Times scaled 840
 \font\sevenrma=\Times scaled 760
 =\Times scaled 680
 \font\fiverma=\Times scaled 600

 \Augment\tenpoint{%
  \textfont0=\tenrma  \scriptfont0=\sevenrma  
  \scriptscriptfont0=\fiverma  }

 \Augment\ninepoint{%
  \textfont0=\ninerma  \scriptfont0=\sevenrma 
  \scriptscriptfont0=\fiverma}

 \Augment\twelvepoint{%
  \textfont0=\twelverma  \scriptfont0=\ninerma  
  \scriptscriptfont0=\sevenrma}

\mathsurround=1pt
\hsize=13.45truecm
\vsize=19.5truecm
\hoffset=1.25truecm
\voffset=2truecm
\advance\baselineskip by 2pt

\predefine\til{\~}
\def\~#1{\relax\ifmmode\widetilde{#1}\else\til{#1}\fi}

\redefine \le{\leqslant}
\redefine \ge{\geqslant}
\define \wt#1{\mathaccent"0365{#1}}
\define \wh#1{\mathaccent"0362{#1}}

\def\U #1{\Underline {#1}}

\define \iss{\,\Mathaccent{\raise -.8 ex\hbox{$\widetilde{}$\kern.1em}}\rightarrow\,}

\define \inlim{{\varinjlim}\vphantom{i}\,}

\define\Car{\mathop{\fam0 C}}

\define \ur{\mathop{\fam0 ur}}

\define \ab{\mathop{\fam0 ab}}

\define \sep{\mathop{\fam0 sep}}

\define \lcf{\mathop{\fam0 lcf}}

\define \kr{\mathop{\fam0 ker}}

\define \chr{\mathop{\fam0 char}\,}
\define \cd{\operatorname{\fam0 cd}}
\define \ind{\operatorname{\fam0 ind}}
\define \pro{\operatorname{\fam0 pro}}

\define \Tr{\operatorname{\fam0 Tr\,}}

\define \Gal{\mathop{\fam0 Gal}}
\define \Hom{\operatorname{\fam0 Hom}}

\define \Spec{\mathop{\fam0 Spec}}

\Mas
\HideEdStuff
\rm 
 

\def\issn{{\nineit ISSN 1464-8997 (on line) 1464-8989 (printed)}}

\def\gtp{{\nineit Published 10 December 2000: \ \copyright\ Geometry \& 
Topology Publications}}

\def\gtv3{{\nineit Geometry \& Topology Monographs, Volume 3 (2000) --
Invitation to higher local fields}}


\def\lione
{{\rms Geometry \& Topology Monographs}}

\def \litwo{{\rms Volume 3: Invitation to higher local fields
}} 

\def\tinfo #1.#2.#3-#4
{{
\noindent  {\lione} \hfill 
\par 
\vskip-1.5pt
\noindent {\litwo} \hfill
\par 
\vskip-1,5pt
\noindent {\rms Part #1, section #2, pages #3--#4} \hfill
\vskip24pt 
}}

\def\tinfos #1.#2.#3-#4
{{
\noindent  {\lione} \hfill 
\par 
\vskip-1.5pt
\noindent {\litwo} \hfill
\par 
\vskip-1.5pt
\noindent {\rms Pages #3--#4} \hfill
\vskip24pt 
}}

\def\tinfoi #1
{{
\noindent  {\lione} \hfill 
\par 
\vskip-1.5pt
\noindent {\litwo} \hfill
\par 
\vskip-1.5pt
\noindent {\rms Pages iii--xi: Introduction and contents} \hfill
\vskip26pt 
}}


  \def\titlepagehead{\hfil}

  \newif\iftitlepage\titlepagefalse
  \newif\ifblankpage\blankpagefalse
  \def\makeheadline{
     \ifblankpage{}\else%
     \iftitlepage
\vbox{\line{\vbox to 8.5pt{}
\ninerm
\copy\HLinebox \hfill
\hglue5mm\ninebf\folio 
\titlepagehead}}%
      \else
\vbox{\ifodd\pageno\rightheadline\else\leftheadline\fi}%
      \fi\vskip 12pt\fi}%
     \def\rightheadline{\line{\vbox to 8.5pt{}%
      \ninerm
\copy\TitleBox \hfill
\hglue5mm\ninebf\folio}}%
     \def\leftheadline{\line{\vbox to 8.5pt{}%
        \unskip\ninerm\unskip\ninebf\folio\hglue5mm
 \hfill \copy\AuthorBox
}}

 \footline={\ifblankpage{}\else
\iftitlepage\ninepoint\sam\hfill
\line{\vbox to 8.5pt{}
\copy\TFLinebox
\hfill
\hglue5mm 
}
            \else
\ninepoint\sam\hfill
\line{\vbox to 8.5pt{}
\copy\FLinebox
\hfill 
\hglue5mm
}
\hfil\fi\global\titlepagefalse\fi}

\def\blankpage{{\blankpagetrue\noindent\hbox to 10pt{\hss}\vfill
\pagebreak}}

\tenpoint\rm 
 